\documentclass[12pt,a4paper]{amsart}
\usepackage{a4}
\usepackage[english]{babel} 
\usepackage{amsmath,amssymb,amsthm}
\def\leq {\leqslant}
\def\le {\leqslant}
\def\ge {\geqslant}
\def\geq {\geqslant}
\def\@bibitem[#1]#2{\item\@biblabel{#1}.\if@filesw
{\def\protect##1{\string##1\space}\immediate\write
\@auxout{\string\bibcite{#2}{#1}}}\fi\ignorespaces\@showtag{#2}}

\theoremstyle{plain}
\newtheorem{theorem}{Theorem}
\newtheorem{rem}{Remark}
\newtheorem{lemma}{Lemma}
\newtheorem{op}{Definition}
\renewcommand{\theequation}%
{\arabic{section}.\arabic{equation}}
\begin{document}

\title{ On embedding theorems of spaces of functions with mixed logarithmic smoothness}
\author{ G. Akishev}
\address{ Lomonosov Moscow University, Kazakhstan Branch \\
Str. Kazhymukan, 11 \\
010010, Astana, Kazakhstan}

\address{
Institute of mathematics and mathematical modeling,\\
Pushkin str, 125, \\
050010, Almaty, \\
 Kazakhstan
 }

\maketitle

\begin{quote}
\noindent{\bf Abstract.}
The article considers the Lorentz space $L_{p,\tau}(\mathbb{T}^{m})$, $2\pi$ of periodic functions of many variables and spaces with mixed logarithmic smoothness. Equivalent norms of a space with mixed logarithmic smoothness are found and embedding theorems are proved.
\end{quote}
\vspace*{0.2 cm}



\section*{Introduction} 

Let $\mathbb{R}^{m}$ ---$m$ be a dimensional Euclidean point space $\overline
{x}=(x_{1},\ldots,x_{m})$ with real coordinates; $\mathbb{T}^{m}=\{\overline x\in \mathbb{R}^{m};\  0\leq x_{j}< 2\pi;\  j=1,\ldots,m\}$  ---$m$ is a dimensional cube and $\mathbb{I}^{m}=[0, 1)^{m}$.

 \smallskip

By $L_{p,\tau}(\mathbb{T}^{m})$ we denote
the Lorentz space of all real-valued Lebesgue measurable functions $f(\overline{x})$ that have $2\pi$-period for each variable and for which the value
\begin{equation*}
\|f\|_{p,\tau} = \left\{\frac{\tau}{p}\int\limits_{0}^{1}\biggl(f^{*}(t)
\biggr)^{\tau}t^{\frac{\tau}{p}-1}dt
\right\}^{\frac{1}{\tau}} , \,\, 1<
p<\infty, 1\leqslant \tau <\infty,
\end{equation*}
is finite, where $f^{*}(y)$ is a non-increasing permutation of the function $|f(2\pi\overline{x})|$,  $\overline{x} \in \mathbb{I}^{m}$ (see \cite[ch.~1, sec. 3, P.~213--216]{1}.

In the case of $\tau=p$, the Lorentz space $L_{p,\tau}(\mathbb{T}^{m})$ coincides with the Lebesgue space $L_{p}(\mathbb{T}^{m})$ with the norm $\|f\|_{p}=\|f\|_{p,p}$
(see \cite[ch.~1, sec.~1.1, p.~11]{2}.

By ${\mathring L}_{p, \tau}(\mathbb{T}^{m})$ we denote the set of all functions $f\in
L_{p, \tau}(\mathbb{T}^{m})$
such that
\begin{equation*}
\int\limits_{0}^{2\pi }f\left( \overline{x} \right) dx_{j}  =0,\;\;
j=1,...,m .
\end{equation*}

We will introduce the following notation:
$a_{\overline{n}}(f)$ -- Fourier coefficients of the function $f\in L_{1}\left(\mathbb{T}^{m}\right)$
according to the system $\{e^{i\langle\overline{n}, 2\pi\overline{x}\rangle}\}_{\mathbb{Z}^{m}},$
where $\mathbb{Z}^{m}$ is a set of points from $\mathbb{R}^{m}$ with integer coordinates.

We will put 
 \begin{equation*}
    \delta _{\overline{s}}(f, 2\pi\overline{x})
=\sum\limits_{\overline{n} \in \rho (\overline{s})
}a_{\overline{n}}(f) e^{i\langle\overline{n} , 2\pi\overline{x}\rangle },
 \end{equation*}
where 
 $\langle\bar{y},\bar{x}\rangle=\sum\limits_{j=1}^{m}y_{j}
x_{j}$, $s_{j} = 1,2,... ,$
 \begin{equation*}
\rho (\overline{s})=\left\{\overline{k} =(k_{1},...,k_{m}) \in \Bbb{Z}^{m}: \, \,  2^{s_{j} -1} \leq \left| k_{j} \right|
<2^{s_{j} } ,j=1,...,m\right\}.
\end{equation*}

Next, $\mathbb{Z}_{+}^{m}$ is a set of points with non-negative integer coordinates.

The value
\begin{equation*}
Y_{l_{1},\ldots,l_{m}}(f)_{p, \tau} = \inf_{T_{l_{j}}} \|f-\sum_{j=1}^{m}T_{l_{j}}\|_{p, \tau}^{*} \;\;,\;\; l_{j} = 0,1,2,...
\end{equation*}
is called the best approximation of the ''angle`` of the function $f\in L_{p, \tau}(\mathbb{T}^{m})$ by trigonometric polynomials, where $T_{l_{j}} \in L_{p}(\mathbb{T}^{m})$a trigonometric polynomial of order $l_{j}$ with respect to the variable $x_{j}, \,\, j = 1,\ldots,m$ (in the case of $\tau=p$, see \cite{3}--\cite{6}).

By $C(p,q,r,y)$ we denote positive values depending on
the parameters indicated in parentheses, generally speaking, different in different
formulas.
For positive values $A(y), B(y)$, write $A\left(y\right) \asymp
B(y)$ means that
there are positive numbers $C_{1},\,C_{2}$ such that $C_{1}\cdot A(y) \leq B(y) \leq C_{2} \cdot A(y)$. For the sake of brevity, in the case of the inequalities $B\ge C_{1}A$ or $B\le C_{2}A$, we will often write $B >> A$ or $B <<A$, respectively.

\begin{op}
Let $k\in \mathbb{N}$ and $\overline{h}=(h_{1}, \ldots , h_{m}) \in\mathbb{R}^{m}$. For the function $f\in L_{1}(\mathbb{T}^{m})$, the total first-order difference at the point $2\pi\overline{x}$ is the value
$\Delta_{\overline{h}}f(2\pi\overline{x}) := f(2\pi\overline{x}+ \overline{h}) - f(2\pi\overline{x})$
and the total difference of order is determined by induction $k$,
\begin{equation*}
\Delta_{\overline{h}}^{k}f(2\pi\overline{x})=\Delta_{\overline{h}}(\Delta_{\overline{h}}^{k - 1}f(2\pi\overline{x})).
\end{equation*}
\end{op}

\begin{op}
The total smoothness modulus of the order $k$ of the function $f\in L_{p, \tau}(\mathbb{T}^{m})$ is the value (in the case of $\tau=p$, see \cite{2}, ch. 4, sec. 4.2)
\begin{equation*}
\omega_{k}(f, t)_{p, \tau} := \sup_{\|\overline{h}\|\leq t}\|\Delta_{\overline{h}}^{k}f\|_{p, \tau}
\end{equation*}
\end{op}

\begin{op}
The mixed modulus of smoothness of a function of order $k$ of the function $f\in L_{p, \tau}(\mathbb{T}^{m})$ is determined by the formula (in the case of $\tau=p$, see \cite{7}, \cite{8}, Ch. 1, sec. 11 )
\begin{equation*}
\omega_{\overline k}(f, \overline{t})_{p, \tau}:=\omega_{k_{1},...k_{m}}(f, t_{1},..., t_{m})_{p, \tau}=\sup_{|h_{1}|\leq t_{1},...,|h_{m}|\leq t_{m}}\|\Delta_{\overline h}^{\overline k} (f)\|_{p, \tau},
\end{equation*}
where 
 $\Delta_{\bar t}^{\bar k} f(2\pi\bar x) = \Delta_{t_{m}}^{k_{m}}(...
\Delta_{t_{1}}^{k_{1}}f(2\pi\overline{x}))$ --- mixed difference of order $\overline{k}$ in increments $\overline{t} = (t_{1},...,t_{m})$.
\end{op}

{\bf The space of the Besov's} (\cite{9}, \cite{10}, \cite{2}, ch. 4, sec. 4. 3).
Let $k \in \mathbb{N}$ and $k > r > 0$, $1\leq \theta \leq \infty$.
The Besov space $\mathbf{B}_{p, \theta}^{r}(\mathbb{T}^{m})$ is called the set of all functions $f\in L_{p}(\mathbb{T}^{m})$, $1\leq p < \infty$ for which
\begin{equation*}
\|f\|_{\mathbf{B}_{p, \theta}^{r}}: = \|f\|_{p} + \Bigl(\int_{0}^{1}(t^{-r}\omega_{k}(f, t)_{p})^{\theta} \frac{dt}{t}\Bigr)^{1/\theta} < \infty, \, \, 
\end{equation*}
for  $1\leq \theta < \infty$ and
\begin{equation*}
\|f\|_{\mathbf{B}_{p, \infty}^{r}}: =\|f\|_{p} +\sup_{t> o}t^{-r}\omega_{k}(f, t)_{p}< \infty
\end{equation*}
for $\theta = \infty$. 

The space $\mathbf{B}_{p, \infty}^{r}(\mathbb{T}^{m})$ is defined and studied by S. M. Nikol'skii \cite{10} and is denoted by the symbol $\mathbf{H}_{p}^{r}(\mathbb{T}^{m})$.

The following generalization of the Besov space is defined in articles \cite{11}--\cite{13}
$\mathbf{B}_{p, \theta}^{r, b}(\mathbb{T}^{m})$ -- the set of all functions $f\in L_{p}(\mathbb{T}^{m})$, $1\leq p < \infty$ for which
\begin{equation*}
\|f\|_{\mathbf{B}_{p, \theta}^{r, b}}: = \|f\|_{p} + \Bigl(\int_{0}^{1}(t^{-r}(1 - \log t)^{b}\omega_{k}(f, t)_{p})^{\theta} \frac{dt}{t}\Bigr)^{1/\theta} < \infty, \, \, 
\end{equation*}
for $0< \theta \leq \infty$, $b> -1/\theta$, $k\in \mathbb{N}$, $k > r>0$.

In the case of $r>0$ and $1<p<\infty$, the norm $\|f\|_{\mathbf{B}_{p, \theta}^{r,b}}$, spaces $\mathbf{B}_{p, \theta}^{r, b}(\mathbb{T}^{m})$
is equivalent to
\begin{equation*}
\|f\|_{p} + \left\{\sum\limits_{s \in \mathbb{Z}_{+}} 2^{sr\theta}(s + 1)^{b\theta} \|\sigma_{s}(f)\|_{p}^{\theta}\right\}^{\frac{1}{\theta}}.
\end{equation*}
 Here  
\begin{equation*}
 \sigma_{s}(f, 2\pi\overline{x}) =  \sum\limits_{\overline{k} \in \square_{2^{s}}\setminus \square_{2^{s-1}}}a_{\overline{k}}(f) e^{i\langle\overline{k} , 2\pi\overline{x}\rangle}, \, \, s\in \mathbb{Z}_{+}, 
\end{equation*}
where 
 \begin{equation*}
 \square_{M} = \{\overline{k} = (k_{1},...,k_{m}) \in \mathbb{Z}^{m} \,\, : |k_{j}| < M, \, \, j=1,...,m\},
 \end{equation*}
for $M\in \mathbb{N}$. 
 
However, in the case of $r=0$, there will be no such equivalence.
Spaces of zero order $B_{p, \theta}^{0}$ from the point of view of approximation theory are defined and investigated by S. M. Nikol'skii \cite{14} and he also noted that there are different definitions of the space $B_{p, \theta}^{0}$. The space $\mathbf{B}_{p, \theta}^{r}(\mathbb{T}^{m})$ for $r=0$ is defined by O. V. Besov \cite{15}.
 
Equivalent norms of the space $\mathbf{B}_{p,\theta}^{0, b}(\mathbb{T}^{m})$
are found in article \cite{13}.

In the works \cite{13}, \cite{16}--\cite{18} the space $B_{p, \theta}^{0, b}(\mathbb{T}^{m})$ is defined as the set of functions $f\in {\mathring L}_{p}(\mathbb{T}^{m})$ for which
\begin{equation*}
\left\{\sum\limits_{s \in \Bbb{Z}_{+}} (s + 1)^{b\theta} \|\sigma_{s}(f)\|_{p, \tau}^{\theta}\right\}^{\frac{1}{\theta}} < \infty,
\end{equation*}
for $1< p < \infty$, $b > -1/\theta$, $0< \theta < \infty$.

Relations between the spaces $B_{p, \theta}^{0, b}(\mathbb{T}^{m})$ and
$\mathbf{B}_{p,\theta}^{0,b}$ investigated in \cite{12}, \cite{16}, \cite{18}, in particular, it is shown that these spaces do not coincide.

{\bf A space with a dominant mixed modulus of smoothness.}
$S_{p}^{\overline{r}}\mathbf{H},$ $S_{p, \theta}^{\overline{r}}\mathbf{B}$ --- the spaces of functions with a dominant mixed derivative are respectively defined by S. M. Nikol'skii \cite{10} and T.I. Amanov (\cite{8} ch. I, p.17). The spaces $S_{p}^{\overline{r}}\mathbf{H},$ $S_{p, \theta}^{\overline{r}}\mathbf{B}$ are called the Nikol'skii--Besov space.

Let $\overline{r} = (r_{1},...,r_{m}),$ $k_{j}>r_{j} > 0$, $j =1,...,m,$
$1 \le p <\infty$, $1\leq \theta \le +\infty.$
The space $S_{p, \theta}^{\overline{r}}\mathbf{B}$
--- consists of all functions $f\in {\mathring L}_{p} (\mathbb{T}^{m})$ for which
\begin{equation*}
\| f\|_{S_{p, \theta}^{\overline{r}}\mathbf{B}
} = \|f\|_{p} + \Biggl[\int_{0}^{1}...\int_{0}^{1}\omega_{\overline{k}}^{\theta}(f, \overline{t})_{p}\prod_{j=1}^{m} \frac{dt_{j}}{t_{j}^{1+\theta r_{j}}}\Biggr]^{\frac{1}{\theta}} < +\infty,
\end{equation*}

P.\,I.~Lizorkin and S.\,M.~Nikol'skii \cite{19} proved that (in the case of $\theta=\infty$, see also \cite{19}, Theorem 4.4.6)
 \begin{equation*}
\| f\|_{S_{p, \theta}^{\overline{r}}\mathbf{B}} \asymp
\left\{\sum\limits_{\overline{s} \in \Bbb{Z}_{+}^{m}} 2^{\langle\overline{s} ,\overline{r}\rangle \theta} \left\| \delta
_{\overline{s}}(f) \right\|_{p }^{\theta}\right\}^{\frac{1}{\theta}}
 \end{equation*}
for  $1 < p < +\infty,$ $1\le \theta \le +\infty.$

M. K. Potapov \cite{4}, \cite{5} defined and investigated the generalization of spaces
$S_{p}^{\overline{r}}\mathbf{H},$ $S_{p,\theta}^{\overline{r}}\mathbf{B}$, with function replacement $t^{r_{j}}$, $j=1,...,m$, to more general functions satisfying certain conditions.

In contrast to \cite{4}, \cite{5}, consider the following space of functions with mixed logarithmic smoothness.

\begin{op} Let $1\leq p< \infty$, $1\leq \tau < \infty$, $0< \theta \leq\infty$, $b_{j}> -1/\theta$, $j=1,...,m$. Through $S_{p, \tau, \theta}^{0, \overline{b}}\mathbf{B}$ denote the space of all functions $f\in{\mathring L}_{p, \tau} (\mathbb{T}^{m})$ for which
\begin{equation*}
\| f\|_{S_{p, \tau, \theta}^{0, \overline{b}}\mathbf{B}}
 := \|f\|_{p, \tau} + \Biggl[\int_{0}^{1}...\int_{0}^{1}\omega_{\overline{k}}^{\theta}(f, \overline{t})_{p, \tau}\prod_{j=1}^{m} \frac{(1-\log t_{j})^{\theta b_{j}}}{t_{j}}dt_{1}...dt_{m}\Biggr]^{\frac{1}{\theta}} < +\infty,
\end{equation*}
where $\overline{b}=(b_{1},...,b_{m})$, $\overline{k}=(k_{1},...,k_{m})$, $k_{j}\in \mathbb{N}$, $j=1,...,m$.  
\end{op}

We note that the space $S_{p, \tau, \theta}^{0, \overline{b}}\mathbf{B}$ is an analog of the space $\mathbf{B}_{p, \tau, \theta}^{0, \overline{b}}$ defined in \cite{21}, using the full smoothness modulus of the function $f\in L_{p, \tau}(\mathbb{T}^{m})$.

We will also consider the following class $S_{p, \tau,\theta}^{0,\overline{b}}B$, consisting of all functions $f\in{\mathring L}_{p, \tau} (\mathbb{T}^{m})$ for which
 \begin{equation*}
\| f\|_{S_{p, \tau, \theta}^{0, \overline{b}}B} =
\left\{\sum\limits_{\bar s \in \Bbb{Z}_{+}^{m}} \prod_{j=1}^{m}(s_{j} + 1)^{b_{j}\theta} \|\delta_{\overline{s}}(f)\|_{p, \tau}^{\theta}\right\}^{\frac{1}{\theta}} < \infty,
 \end{equation*}
 for $1 < p < +\infty,$ $1\leq \tau< \infty$, $1\le \theta \le +\infty$, $\overline{b}=(b_{1},...,b_{m})$, $b_{j}b_{j}> -1/\theta$, $j=1,...,m$.

In the case of $\tau=p$, in the place of $S_{p, p,\theta}^{0, \overline{b}}B$, $S_{p, p, \theta}^{0, \overline{b}}\mathbf{B}$ we will write $S_{p, \theta}^{0, \overline{b}}B$, $S_{p, \theta}^{0, \overline{b}}\mathbf{B}$ respectively.

We note that the classes $S_{p, \theta}^{0, \overline{b}}\mathbf{B}$ and $S_{p,\theta}^{0, \overline{b}}B$ are analogs of the spaces $\mathbf{B}_{p, \tau, \theta}^{0, \overline{b}}$ and $B_{p, \tau, \theta}^{0, \overline{b}}$ , defined and investigated in \cite{12}, \cite{13}, \cite{17}, \cite{18}.

The main purpose of the article is to find equivalent norms of the space $S_{p, \theta}^{0, \overline{b}}\mathbf{B}$ and prove embedding theorems for the spaces $S_{p, \theta}^{0, \overline{b}}\mathbf{B}$ and $S_{p, \theta}^{0, \overline{b}}B$.

The article consists, in addition to the introduction, of two sections. In the first section, the auxiliary statements necessary to prove the main results are proved.
The main results are formulated and proved in the second section.
Theorems 1, 2, 4 of this section are analogs, respectively, of formula 3.7 in \cite[p. 1043]{21}, Theorem 4.3 in \cite{13}, Theorem 4.1 in \cite{16} and Theorem 3.3 in \cite{12}. The 3rd point of Theorem 5 corresponds to corollary 3.6 in \cite[p. 1047]{21}.

  \smallskip
\setcounter{equation}{0}
\setcounter{lemma}{0}
\setcounter{theorem}{0}

\section{Auxiliary statements} 

\smallskip

First, we introduce additional notation and give auxiliary statements. In the future, we will use the following notation.

The set of indices $\{1,...,m\}$ is denoted by the symbol $e_{m},$ its arbitrary subset - through $e$ and $|e|$ - the number of elements of $e$ .

If given an element $\overline{r} =(r_{1},\ldots,r_{m})$, $m$ -- dimensional space with
non-negative coordinates, then $\overline{r}^{e} =(r_{1}^{e},\ldots,r_{m}^{e})$
vector with components $r_{j}^{e} = r_{j}$ for $j\in e$ and $r_{j}^{e}= 0$ for $j\notin e.$

Let $\overline{l} = (l_{1},\ldots,l_{m})$ be an element of $m$ --dimensional space with
positive integer coordinates and a nonempty set $e\subset e_{m}.$
We set 
 $$
G_{\overline{l}}(e) = \{ \overline{k} = (k_{1},\ldots,k_{m})\in \mathbb{Z}^{m} :
|k_{j}|\le l_{j}, j\in e \quad |k_{j}| > l_{j}, j\notin e \}.
 $$
For the given numbers $b_{\overline{n}}$, the mixed difference is determined
by the formula
 $$
\Delta {b_{\overline{n}}} = \sum\limits_{\overline{0} \le \overline\varepsilon \le \overline{1}}(-1)^{m - \sum\limits_{j=1}^{m} \varepsilon_{j}} b_{\overline{n} - \overline{1} + \overline\varepsilon},
 $$
where 
 $\overline\varepsilon = (\varepsilon_{1},\ldots,\varepsilon_{m})$ and
$\varepsilon_{j}=0$ or $\varepsilon_{j}=1,$ 
$\overline{n} - \overline{1} + \overline\varepsilon = (n_{1} - 1 + \varepsilon_{1},\ldots,n_{m} - 1 + \varepsilon_{m}).$

We consider partial sums over various variables:
 $$
S_{\overline{l}}(f, 2\pi\overline{x}) = S_{l_{1},\ldots,l_{m}}(f, 2\pi\overline{x}) =
\sum\limits_{|k_{1}|\le l_{1}}\ldots\sum\limits_{|k_{m}|\le l_{m} }
a_{\overline{k}} (f) e^{i\langle\overline{k} , 2\pi\overline{x}\rangle}
 $$
-- partial sum for all variables;
 $$
S_{l_{1}, \infty}(f, 2\pi\overline{x}) =
\sum\limits_{|k_{1}|\le l_{1}}\sum\limits_{k_{2}=-\infty }^{+\infty}\ldots\sum\limits_{k_{m}=-\infty}^{+\infty}
a_{\overline{k}}(f) e^{i\langle\overline{k} , 2\pi\overline{x}\rangle}
 $$
-- partial sum of the variable $ x_{1}\in [0, 1).$ 
In the more general case 
 $$
S_{\overline{{l}^{e}, \infty}}(f, 2\pi\overline{x}) =
\sum\limits_{\overline{k} \in \prod_{j \in e} [-l_{j}, l_{j}]\times \mathbb{R}^{m-|e|} }
a_{\overline{k}}(f) e^{i\langle\overline{k} , 2\pi\overline{x}\rangle}
 $$
-- partial sum of variables $x_{j}\in [0, 1)$ for $j\in e.$ 

For a given subset of $e\subset e_{m}$, we put 
 $$
U_{\overline{l}}(f, 2\pi\overline{x}) = \sum\limits_{e\subset e_{m}, e \neq \emptyset} \;\; \sum\limits_{\overline{k} \in G_{\overline{l}}(e)} a_{\overline{k} } (f) e^{i\langle\overline{k} ,2\pi\overline{x}\rangle}.
 $$
In particular, for $m=2$ we have (see for example \cite{7} )
 $$
U_{l_{1},l_{2}}(f, 2\pi\overline{x}) = S_{l_{1} ,\infty }(f, 2\pi\overline{x}) + S_{\infty,l_{2}  }(f, 2\pi\overline{x}) - S_{l_{1},l_{2}  }(f, 2\pi\overline{x}).
 $$
We give some properties of the mixed modulus of smoothness of the function.
\begin{lemma}\label{lem1} 
Let $1 < p, \tau < +\infty,$ $\alpha_{j}\in \mathbb{N}$, for $j=1, \ldots , m$ and the function $f, g\in L_{p, \tau}(\mathbb{T}^{m})$. Then

1. $\omega_{\overline{\alpha}}(f + g, \delta_{1}, ...,\delta_{m})_{p, \tau}<<\omega_{\overline{\alpha}}(f, \delta_{1}, ...,\delta_{m})_{p, \tau} + \omega_{\overline{\alpha}}(g, \delta_{1}, ...,\delta_{m})_{p, \tau}$;
 
2. $\omega_{\overline{\alpha}}(f, \delta_{1}, ...,\delta_{m})_{p, \tau}$ does not decrease for each variable $\delta_{j}\geq 0$, $j=1,...,m$; 

3. $\omega_{\overline{\alpha}}(f, \lambda_{1}\delta_{1}, ...,\lambda_{m}\delta_{m})_{p, \tau} << \prod_{j=1}^{m}\lambda_{j}^{\alpha_{j}}\omega_{\overline{\alpha}}(f, \delta_{1}, ...,\delta_{m})_{p, \tau}$,
for numbers $\lambda_{j}\geq 0$, $j=1,...,m$;

4. for a trigonometric polynomial  
 \begin{equation*}
T_{\overline{n}}(2\pi\overline{x}) = \sum\limits_{k_{1}=-n_{1}}^{n_{1}} \ldots \sum\limits_{k_{m}=-n_{m}}^{n_{m}}c_{\overline{k}}e^{i\langle\overline{k} , 2\pi\overline{x}\rangle}, \, \, n_{j}\in \mathbb{Z}_{+}, j=1,...,m, \, \, \overline{x}\in \mathbb{I}^{m} 
\end{equation*} 
and its derivative $T_{\overline{n}}^{(\alpha_{1},...,\alpha_{m})}(2\pi\overline{x})$, the inequality holds
 \begin{equation*}
 \omega_{\overline{\alpha}}(T_{\overline{n}}, \delta_{1}, ...,\delta_{m})_{p, \tau}<< 
\prod_{j=1}^{m}\delta_{j}^{\alpha_{j}}\|T_{\overline{n}}^{(\alpha_{1},...,\alpha_{m})}\|_{p, \tau}. 
 \end{equation*}
 \end{lemma}
\proof This lemma is proved as Theorem 4.1 in \cite{7}. To prove property 4, we use the well-known equality (see for example \cite{8}, Ch. 1, sec. 11)
\begin{equation*}
\Delta_{\overline{h}}^{\overline{\alpha}}T_{\overline{n}}(2\pi\overline{x}) = \prod_{i=1}^{m}\int_{0}^{h_{i}} \ldots \int_{0}^{h_{i}}T_{\overline{n}}^{(\alpha_{1},...,\alpha_{m})}(2\pi x_{1}+\sum\limits_{j=1}^{\alpha_{1}}u_{j}^{1}, \ldots 2\pi x_{1}+\sum\limits_{j=1}^{\alpha_{m}}u_{j}^{m})du_{1}^{i} \ldots du_{\alpha_{i}}^{i} 
\end{equation*}  
and the generalized Minkowski's inequality in the Lorentz space . Then
\begin{equation*}
\|\Delta_{\overline{h}}^{\overline{\alpha}}T_{\overline{n}}\|_{p, \tau}<< 
\prod_{i=1}^{m}\Bigl |\int_{0}^{h_{i}} \ldots \int_{0}^{h_{i}}\|T_{\overline{n}}^{(\alpha_{1},...,\alpha_{m})}\|_{p, \tau} du_{1}^{i} \ldots du_{\alpha_{i}}^{i}\Bigr |= C\prod_{i=1}^{m}|h_{i}|^{\alpha_{i}}\|T_{\overline{n}}^{(\alpha_{1},...,\alpha_{m})}\|_{p, \tau},
\end{equation*}
for 
$h_{i}\in \mathbb{R}$, $i=1,...,m$. Hence the inequality in property 4.  \hfill $\Box$

\begin{rem} In the case of $\tau =p$, Lemma 1 was previously proved in \cite[theorem 5.1]{7}.
\end{rem}

\begin{lemma}\label{lem2} (Bernstein's Inequality).
Let $1 <p, \to < +\infty,$ $\alpha_{j}\in \mathbb{Z}_{+}$, for $j=1, \ldots , m$. Then for the trigonometric polynomial $T_{\overline{n}}$ there is an inequality
\begin{equation*}
\|T_{\overline{n}}^{(\alpha_{1},...,\alpha_{m})}\|_{p, \tau}<<\prod_{j=1}^{m}(n_{j}+1)^{\alpha_{j}} \|T_{\overline{n}}\|_{p, \tau}
\end{equation*}
\end{lemma}
\proof 
In the Lebesgue space $L_{p}(\mathbb{T}^{m})$, $1\leq p\leq\infty$, Bernstein's inequality is known (\cite{2}, ch. 2, sec. 2.7. 2)
\begin{equation*}
\|T_{\overline{n}}^{(\alpha_{1},...,\alpha_{m})}\|_{p}<<\prod_{j=1}^{m}(n_{j}+1)^{\alpha_{j}} \|T_{\overline{n}}\|_{p}.
\end{equation*}
It is also known that the Lorentz space $L_{p,\tau}(\mathbb{T}^{m})$ is an interpolation space between $L_{p_{1}}(\mathbb{T}^{m})$ and $L_{p_{2}}(\mathbb{T}^{m})$, $1< p_{1}<p <p_{2}< \infty$ (for example, see \cite[p. 132]{20}). Therefore, from the previous inequality we obtain the statement of the lemma.  \hfill $\Box$

\begin{lemma}\label{lem3} 
Let $1 < p < +\infty,$\,\,$1 < \tau < +\infty$ and $f\in {\mathring L}_{p, \tau}(\mathbb{T}^{m})$.  Then 
 \begin{equation*}
\|f - U_{l_{1},...,l_{m}} (f)\|_{p, \tau}\ll Y_{l_{1},...,l_{m}}(f)_{p, \tau}.
\end{equation*}
\end{lemma}
\proof The proof is carried out according to the scheme \cite{3}. Let $m=2$ and $T_{l_{1}}\in
{\mathring L}_{p, \tau}(\mathbb{T}^{2}) $ is a trigonometric polynomial of order
$l_{1}$ by $x_{1}$, $T_{l_{2}} \in {\mathring L}_{p, \tau}(\mathbb{T}^{m})$--
a trigonometric polynomial of order $l_{2}$ by $x_{2}$ and
$T_{l_{1},l_{2}} \in {\mathring L}_{p, \tau}(\mathbb{T}^{m})$ -- a trigonometric
polynomial of order $l_{1}$ by $x_{1}$ and $l_{2}$ by $x_{2}.$ It is known that
 \begin{equation*}
S_{l_{1} ,\infty }(T_{l_{1}}, 2\pi\bar x) = T_{l_{1}}(2\pi\bar x),\;\; S_{\infty,l_{2}  }(T_{l_{2}}, 2\pi\bar x) = T_{l_{2}}(2\pi\bar x),
\end{equation*}
 \begin{equation*}
 S_{l_{1} ,\infty }(T_{l_{2}}, 2\pi\bar x) = S_{l_{1},l_{2}  }(T_{l_{2}}, 2\pi\bar x),\;\;
S_{\infty,l_{2}  }(T_{l_{1}}, 2\pi\bar x) = S_{l_{1},l_{2}  } (T_{l_{1}}, 2\pi\bar x),
\end{equation*}
 \begin{equation*}
T_{l_{1},l_{2}}(2\pi\bar x) = S_{l_{1},l_{2}  }( T_{l_{1},l_{2}}, 2\pi\bar x) =
S_{l_{1},\infty }(T_{l_{1},l_{2}}, 2\pi\bar x)=S_{\infty,l_{2}  }(T_{l_{1},l_{2}}, 2\pi\bar x), \, \, \overline{x} \in \mathbb{I}^{m}.
 \end{equation*}
Then  
$f(2\pi\overline{x}) - U_{l_{1},...,l_{m}} (f, 2\pi\overline{x} ) = \varphi (2\pi\overline{x}) - S_{l_{1},\infty }(\varphi, 2\pi\overline{x})
- S_{\infty,l_{2}}(\varphi, 2\pi\overline{x}) + S_{l_{1},l_{2}}(\varphi, 2\pi\overline{x}),$ where 
$\varphi (2\pi\overline{x}) = f(2\pi\overline{x}) - T_{l_{1},\infty}(2\pi\overline{x}) - T_{\infty,l_{2}}(2\pi\overline{x}) + T_{l_{1},l_{2}}(2\pi\overline{x})$, $\overline{x} \in \mathbb{I}^{m}$. 
According to the boundedness of the rectangular partial sum operator in space
Lorentz we have
 \begin{equation*}
\|S_{l_{1} ,\infty }(\varphi)\|_{p, \tau} \ll
\|\varphi\|_{p, \tau} , \;\;
\|S_{\infty ,l_{2} }(\varphi)\|_{p, \tau} \ll \|\varphi\|_{p, \tau},
 \end{equation*}
 \begin{equation*}
\|S_{l_{1},l_{2} }(\varphi)\|_{p, \tau} \ll
\|\varphi\|_{p, \tau}, \; \; 1<p, \tau < \infty.
 \end{equation*}
Therefore , by the metric property we get 
$\|f - U_{l_{1},l_{2}} (f )\|_{p, \tau}\ll \|\varphi\|_{p, \tau}.$
Since $T_{l_{1},\infty}, T_{\infty,l_{2}}, T_{l_{1},l_{2}}$ are any , then the statement of the lemma follows.  \hfill $\Box$

\begin{lemma}\label{4} (Direct theorem). 
If $f\in \mathring{L}_{p, \tau}(\mathbb{T}^{m})$,  $1<p<+\infty,$
$1<\tau <+\infty$, then
 \begin{equation*}
Y_{\overline{n}}(f)_{p, \tau}\ll \omega_{\bar k}
\left(f,\frac{1}{n_{1}+1},...,\frac{1}{n_{m}+1}\right)_{p, \tau}.
 \end{equation*}
\end{lemma}
 \proof  We will prove the lemma as in the case of the space $L_{p}(\mathbb{T}^{m})$ in [3]. We will consider the  well-known Jackson--Stechkin kernel (see \cite{23}, sec. 1, subsec. 2) 
  \begin{equation*}
F_{l}(u)=b_{r}\cdot \left(\frac{\sin \frac{ru}{2}}{\sin \frac{u}{2}}
\right)^{2k_{0}},
 \end{equation*}
where the natural number $k_{0}>\frac{k+1}{2}$, and the natural number $r$ is chosen so
that $\frac{l}{2k_{0}}<r\leq\frac{l}{2k_{0}}+1$ and
 \begin{equation*}
b_{r}=\left[\int_{-\pi}^{\pi}\left(\frac{\sin \frac{ru}{2}}{\sin \frac{u}{2}}
\right)^{2k_{0}}du\right]^{-1},
\end{equation*}
so 
 $\int_{-\pi}^{\pi}F_{l}(u)du=1.$
It is known that for $\mu\in[0,k]$ the inequality holds
 \begin{equation*}
\int_{-\pi}^{\pi}F_{l}(u)|u|^{\mu}du\ll (l+1)^{-\mu}.
\end{equation*}
Now we will consider the function

$$
T_{l_{j}}(f, 2\pi\overline{x})
$$
$$=(-1)^{k_{j}+1}\int_{-\pi}^{\pi}F_{l_{m}}(t)
\sum_{\nu_{j}=1}^{k_{j}}(-1)^{k_{j}-\nu_{j}} C_{k_{j}}^{\nu_{j}}
  f(2\pi x_{1},...,2\pi x_{j-1}, 2\pi x_{j}+t_{j}\nu_{j}, 2\pi x_{j+1},...,2\pi x_{m})dt_{j}.
$$
It is known that $T_{l_{j}}(f, 2\pi\overline{x})$ is a measurable function of $m$ variables,
$2\pi$ is periodic for each variable and is a trigonometric
polynomial of the order $l_{j}$ for the variable $2\pi x_{j}$ (see \cite{23}, sec. 1, subsec. 2).
We now take the polynomial $T_{l_{m}}(f, 2\pi\overline{x})$ and form the function
$\varphi_{m}(2\pi \overline{x})=f(2\pi\overline{x})-T_{l_{m}}(f, 2\pi\overline{x}).$
For this function, we construct a polynomial $T_{l_{m-1}}(\varphi_{m}(2\pi\overline{x}))$ and consider the function
 $ 
\varphi_{m-1,m}(2\pi\overline{x}) = \varphi_{m}(2\pi\overline{x}) -
T_{l_{m-1}}(\varphi_{m}(2\pi\overline{x})).
 $ 
Next, we make a polynomial $T_{l_{m-2}}(\varphi_{m-1,m}, 2\pi\overline{x})$ and consider
the function
\begin{equation*}
\varphi_{m-2,m-1,m}(2\pi\overline{x})=\varphi_{m-1,m}(2\pi\overline{x})-
T_{l_{m-2}}(\varphi_{m-1,m},2\pi \overline{x}).
 \end{equation*}
Continuing this process we will get the function
 $$
\varphi_{1,2,...,m}(2\pi\overline{x})=\varphi_{2,...,m}(2\pi\overline{x})-
T_{l_{1}}(\varphi_{2,...,m}, 2\pi\overline{x}).
 $$
Then 
 $
\varphi_{1,...,m}(2\pi\overline{x})=f(2\pi\overline{x})-
\left[T_{l_{1}}(\varphi_{2,...,m}, 2\pi\overline{x})+...+T_{l_{m}}(f, 2\pi\overline{x})\right].
 $
Therefore, 
\begin{equation}\label{eq1 1}
Y_{l_{1},...,l_{m}}(f)_{p, \tau}\leq \left\|f-
\left[T_{l_{1}}(\varphi_{2,...,m}, \overline{x})+...
+T_{l_{m}}(f, \overline{x})\right]\right\|_{p, \tau}. 
 \end{equation}
Taking into account the properties of the kernels $F_{l_{j}}(t)$ and the way the polynomials are formed $T_{l_{j}}$
you can make sure that
\begin{equation*}
f(2\pi\overline{x})-\sum_{j=1}^{m}T_{l_{j}}(\varphi_{j+1,...,m}, 2\pi\overline{x})=\int_{-\pi}^{\pi}...\int_{-\pi}^{\pi}\prod_{j=1}^{m}F_{l_{j}}(t_{j})
\Delta_{t_{1},...,t_{m}}^{k_{1},...,k_{m}}(f, 2\pi\overline{x})dt_{1}...dt_{m}.
 \end{equation*}
Therefore , from (1. 1) by virtue of the generalized Minkowski inequality in space
We will get Lorentz
 \begin{multline}\label{eq1 2}
Y_{l_{1},...,l_{m}}(f)_{p, \tau}\leq \left\|
\int_{-\pi}^{\pi}...
\int_{-\pi}^{\pi}\prod_{j=1}^{m}F_{l_{j}}(t_{j})
\Delta_{t_{1},...,t_{m}}^{k_{1},...,k_{m}}(f, 2\pi\overline{x})dt_{1}...dt_{m}
\right\|_{p, \tau}
 \\
\leq \int_{-\pi}^{\pi}...
\int_{-\pi}^{\pi}\prod_{j=1}^{m}F_{l_{j}}(t_{j})
\omega_{k_{1},...,k_{m}}(f, |t_{1}|,...,|t_{m}|)_{p, \tau}
dt_{1}...dt_{m}. 
 \end{multline}
Now by the property of the smoothness module and kernels $F_{l_{j}}(t_{j})$ we have
  \begin{multline*}
\int_{-\pi}^{\pi}...
\int_{-\pi}^{\pi}\prod_{j=1}^{m}F_{l_{j}}(t_{j})
\omega_{k_{1},...,k_{m}}(f,|t_{1}|,...,|t_{m}|)_{p, \tau}
dt_{1}...dt_{m}
 \\
=\int_{-\pi}^{\pi}...
\int_{-\pi}^{\pi}\prod_{j=1}^{m}F_{l_{j}}(t_{j})\cdot
\omega_{k_{1},...,k_{m}}(f, \frac{(l_{1}+1)|t_{1}|}{l_{1}+1},...,
\frac{(l_{m}+1)|t_{m}|}{l_{m}+1})_{p, \tau}
dt_{1}...dt_{m} 
\\
\leq \int_{-\pi}^{\pi}...
\int_{-\pi}^{\pi}\prod_{j=1}^{m}F_{l_{j}}(t_{j})
\left(|t_{j}|(l_{j}+1)+1\right)^{k_{j}}
dt_{1}...dt_{m}
\\
\times \omega_{k_{1},...,k_{m}}(f,\frac{1}{l_{1}+1},...,
\frac{1}{l_{m}+1})_{p, \tau} 
\ll \omega_{k_{1},...,k_{m}}(f,\frac{1}{l_{1}+1},...,
\frac{1}{l_{m}+1})_{p, \tau}.
 \end{multline*}
 \hfill $\Box$

\begin{lemma}\label{lem5} (Inverse theorem). 
If $f\in {\mathring L}_{p, \tau}(\mathbb{T}^{m})$,  $1<p<+\infty,$
$1<\tau <+\infty$, $\alpha_{j}\in \mathbb{N}$, then
 \begin{equation*}
 \omega_{\bar\alpha}\bigl(f,\frac{1}{n_{1}+1},...,\frac{1}{n_{m}+1}\bigr)_{p, \tau}
\ll\prod_{j=1}^{m}n_{j}^{-\alpha_{j}}\sum\limits_{\nu_{1}=1}^{n_{1}+1} \ldots \sum\limits_{\nu_{m}=1}^{n_{m}+1} \prod_{j=1}^{m}\nu_{j}^{\alpha_{j} - 1}Y_{\overline{\nu}}(f)_{p, \tau}.
 \end{equation*}
\end{lemma}
\proof The lemma is proved according to a well-known scheme (see for example \cite{5}, \cite{7}). Therefore, we will prove the lemma for $m=2$.
Let
 $f\in {\mathring L}_{p, \tau}(\mathbb{T}^{m})$.
 For $n_{j}\in \mathbb{N}$, we choose a non-negative integer $k_{j}$ such that $2^{k_{j}}\leq n_{j}< 2^{k_{j} + 1}$, $j=1, \ldots , m$.
Then by the properties of the mixed modulus of smoothness we get 
 \begin{multline}\label{eq1 3}
 \omega_{\bar\alpha}\bigl(f,\frac{1}{n_{1}+1}, \frac{1}{n_{2}+1}\bigr)_{p, \tau}\leq 
 \omega_{\bar\alpha}\bigl(f,\frac{1}{2^{k_{1}}}, \frac{1}{2^{k_{2}}}\bigr)_{p, \tau}
\\
\leq \omega_{\bar\alpha}\bigl(f-S_{2^{k_{1}}, \infty}(f)-S_{\infty, 2^{k_{2}}}(f)+S_{2^{k_{1}}, 2^{k_{2}}}(f), \frac{1}{2^{k_{1}}}, \frac{1}{2^{k_{2}}}\bigr)_{p, \tau}
\\
+ \omega_{\bar\alpha}\bigl(S_{2^{k_{1}}, \infty}(f)-S_{\infty, 2^{k_{2}}}(f), \frac{1}{2^{k_{1}}}, \frac{1}{2^{k_{2}}}\bigr)_{p, \tau} + \omega_{\bar\alpha}\bigl(S_{2^{k_{1}}, 2^{k_{2}}}(f), \frac{1}{2^{k_{1}}}, \frac{1}{2^{k_{2}}}\bigr)_{p, \tau}
\\
\leq 2^{\alpha_{1} + \alpha_{2}}\|f-S_{2^{k_{1}}, \infty}(f)-S_{\infty, 2^{k_{2}}}(f)+S_{2^{k_{1}}, 2^{k_{2}}}(f)\|_{p, \tau} + \omega_{\bar\alpha}\bigl(S_{2^{k_{1}}, \infty}(f-S_{\infty, 2^{k_{2}}}(f)), \frac{1}{2^{k_{1}}}, \frac{1}{2^{k_{2}}}\bigr)_{p, \tau}
 \\
+ \omega_{\bar\alpha}\bigl(S_{\infty, 2^{k_{2}}}(f-S_{2^{k_{1}}, \infty}(f)), \frac{1}{2^{k_{1}}}, \frac{1}{2^{k_{2}}}\bigr)_{p, \tau} + \omega_{\bar\alpha}\bigl(S_{2^{k_{1}}, 2^{k_{2}}}(f), \frac{1}{2^{k_{1}}}, \frac{1}{2^{k_{2}}}\bigr)_{p, \tau}
\\
\ll \|f-S_{2^{k_{1}}, \infty}(f)-S_{\infty, 2^{k_{2}}}(f)+S_{2^{k_{1}}, 2^{k_{2}}}(f)\|_{p, \tau} + 2^{-k_{1}\alpha_{1}}\|S^{(\alpha_{1}, 0)}_{2^{k_{1}}, \infty}(f-S_{\infty, 2^{k_{2}}}(f))\|_{p, \tau}
\\
+2^{-k_{2}\alpha_{2}}\|S^{(0, \alpha_{2})}_{\infty, 2^{k_{2}}}(f-S_{2^{k_{1}}, \infty}(f))\|_{p, \tau} + 
2^{-(k_{1}\alpha_{1}+k_{2}\alpha_{2})}\|S^{(\alpha_{1}, \alpha_{2})}_{2^{k_{1}}, 2^{k_{2}}}(f)  \|_{p, \tau}.
  \end{multline}
Further, by the norm property, Bernstein's inequality (Lemma 2) and by Lemma 3 we will have
\begin{multline*}
2^{-(k_{1}\alpha_{1}+k_{2}\alpha_{2})}\|S^{(\alpha_{1}, \alpha_{2})}_{2^{k_{1}}, 2^{k_{2}}}(f)  \|_{p, \tau}\ll 2^{-(k_{1}\alpha_{1}+k_{2}\alpha_{2})} \sum\limits_{\mu_{1}=0}^{k_{1}}\sum\limits_{\mu_{2}=0}^{k_{2}}\|\delta_{\mu_{1}, \mu_{2}}  ^{(\alpha_{1}, \alpha_{2})}(f)\|_{p, \tau}  \ll 2^{-(k_{1}\alpha_{1}+k_{2}\alpha_{2})}
\\
\times
\sum\limits_{\mu_{1}=0}^{k_{1}}\sum\limits_{\mu_{2}=0}^{k_{2}}2^{\mu_{1}\alpha_{1}+\mu_{2}\alpha_{2}}\|\delta_{\mu_{1}, \mu_{2}}  (f)\|_{p, \tau}\ll 2^{-(k_{1}\alpha_{1}+k_{2}\alpha_{2})}
\end{multline*}
 \begin{multline}\label{eq1 4}
\times
\sum\limits_{\mu_{1}=0}^{k_{1}}\sum\limits_{\mu_{2}=0}^{k_{2}}2^{\mu_{1}\alpha_{1}+\mu_{2}\alpha_{2}}\|f - U_{[2^{\mu_{1} - 1}], [2^{\mu_{2} - 1}]}(f) \|_{p, \tau} 
\\
\ll 2^{-(k_{1}\alpha_{1}+k_{2}\alpha_{2})}\sum\limits_{\mu_{1}=0}^{k_{1}}\sum\limits_{\mu_{2}=0}^{k_{2}}2^{\mu_{1}\alpha_{1}+\mu_{2}\alpha_{2}}Y_{[2^{\mu_{1} - 1}], [2^{\mu_{2} - 1}]}(f)_{p, \tau}
\end{multline}
Using the norm property, the Bernstein inequality and the Littlewood--Paley theorem in the Lorentz space \cite{24} we obtain
\begin{multline*} 
2^{-k_{1}\alpha_{1}}\|S^{(\alpha_{1}, 0)}_{2^{k_{1}}, \infty}(f-S_{\infty, 2^{k_{2}}}(f))\|_{p, \tau} =2^{-k_{1}\alpha_{1}}\Bigl\|\sum\limits_{\mu_{1}=0}^{k_{1}}\delta_{\mu_{1}}  ^{(\alpha_{1}, 0)}(f - S_{\infty, 2^{k_{2}}}(f))\Bigr\|_{p, \tau}
\\
\ll 2^{-k_{1}\alpha_{1}}\sum\limits_{\mu_{1}=0}^{k_{1}}2^{\mu_{1}\alpha_{1}}\|\delta_{\mu_{1}}(f - S_{\infty, 2^{k_{2}}}(f))\|_{p, \tau} = C2^{-k_{1}\alpha_{1}}\sum\limits_{\mu_{1}=0}^{k_{1}}2^{\mu_{1}\alpha_{1}}\Bigl\|\sum\limits_{\mu_{2}=k_{2}+1}^{\infty}\delta_{\mu_{1}, \mu_{2}}(f)\Bigr\|_{p, \tau}
\\
\ll 2^{-k_{1}\alpha_{1}}\sum\limits_{\mu_{1}=0}^{k_{1}}2^{\mu_{1}\alpha_{1}}\Bigl\|\Bigl(\sum\limits_{\mu_{2}=k_{2}+1}^{\infty}|\delta_{\mu_{1}, \mu_{2}}(f)|^{2}\Bigr)^{1/2}\Bigr\|_{p, \tau}
\end{multline*}
\begin{multline} \label{eq1 5}
 \ll 2^{-k_{1}\alpha_{1}}\sum\limits_{\mu_{1}=0}^{k_{1}}2^{\mu_{1}\alpha_{1}}\Bigl\|\Bigl(\sum\limits_{\mu_{2}=k_{2}+1}^{\infty}\sum\limits_{\nu_{1}=\mu_{1}}^{\infty}|\delta_{\nu_{1}, \mu_{2}}(f)|^{2}\Bigr)^{1/2}\Bigr\|_{p, \tau} \ll
 2^{-k_{1}\alpha_{1}}
\\
\times 
 \sum\limits_{\mu_{1}=0}^{k_{1}+1}2^{\mu_{1}\alpha_{1}}Y_{[2^{\mu_{1} - 1}], [2^{k_{2}}]}(f)_{p, \tau} \ll 2^{-(k_{1}\alpha_{1}+k_{2}\alpha_{2})}\sum\limits_{\mu_{1}=0}^{k_{1}+1}\sum\limits_{\mu_{2}=0}^{k_{2} + 1}2^{\mu_{1}\alpha_{1}+\mu_{2}\alpha_{2}}Y_{[2^{\mu_{1} - 1}], [2^{\mu_{2} - 1}]}(f)_{p, \tau}.
\end{multline}
Similarly, it is proved that
 \begin{equation}\label{eq1 6}
2^{-k_{2}\alpha_{2}}\|S^{(0, \alpha_{2})}_{\infty, 2^{k_{2}}}(f-S_{2^{k_{1}}, \infty}(f))\|_{p, \tau} \ll 2^{-(k_{1}\alpha_{1}+k_{2}\alpha_{2})}\sum\limits_{\mu_{1}=0}^{k_{1}+1}\sum\limits_{\mu_{2}=0}^{k_{2} + 1}2^{\mu_{1}\alpha_{1}+\mu_{2}\alpha_{2}}Y_{[2^{\mu_{1} - 1}], [2^{\mu_{2} - 1}]}(f)_{p, \tau}.
\end{equation}
Further on Lemma 3 we have
 \begin{equation}\label{eq1 7}
\|f-S_{2^{k_{1}}, \infty}(f)-S_{\infty, 2^{k_{2}}}(f)+S_{2^{k_{1}}, 2^{k_{2}}}(f)\|_{p, \tau} \ll Y_{[2^{k_{1}}], [2^{k_{2}}]}(f)_{p, \tau}.
\end{equation}
Now from the inequalities (1. 3)--(1. 7) we get
 \begin{multline*}
  \omega_{\bar\alpha}\bigl(f,\frac{1}{n_{1}+1}, \frac{1}{n_{2}+1}\bigr)_{p, \tau}
 \ll 2^{-(k_{1}\alpha_{1}+k_{2}\alpha_{2})}\sum\limits_{\mu_{1}=0}^{k_{1}+1}\sum\limits_{\mu_{2}=0}^{k_{2} + 1}2^{\mu_{1}\alpha_{1}+\mu_{2}\alpha_{2}}Y_{[2^{\mu_{1} - 1}], [2^{\mu_{2} - 1}]}(f)_{p, \tau}
 \\
 \ll n_{1}^{-\alpha_{1}}n_{2}^{-\alpha_{2}}\sum\limits_{\nu_{1}=1}^{n_{1}+1}\sum\limits_{\nu_{2}=1}^{n_{2} + 1}\nu_{1}^{\alpha_{1} - 1}\nu_{2}^{\alpha_{2}-1}Y_{\nu_{1} - 1, \nu_{2} - 1}(f)_{p, \tau}.
\end{multline*}
 \hfill $\Box$

{\bf Remark 2.} In the case of $\tau = p$, Lemmas 3--5 are proved in \cite{3}, \cite{7}, and Lemma 1 --- in \cite{25}. 
  \smallskip

\setcounter{equation}{0}
\setcounter{lemma}{0}
\setcounter{theorem}{0}

\section{Embedding theorems of classes with mixed logarithmic smoothness} 

\smallskip
  
  \begin{theorem}\label{th1}
  Let $1 <p <+ \infty$, $1< \tau < \infty$, $0< \theta \leq \infty$ and the numbers $b_{j}> - \frac{1}{\theta}$, for $j=1,\ldots, m$. Then for the function $f\in S_{p, \tau, \theta}^{\bar r}\mathbf{B}$ the relation 
     \begin{equation}\label{eq2 1}
 \| f\|_{S_{p, \tau, \theta}^{0, \bar b}\mathbf{B}}
\asymp \Biggl(\sum\limits_{\nu_{m} =0}^{\infty}...\sum\limits_{\nu_{1} =0}^{\infty}\prod_{j=1}^{m}(\nu_{j} + 1)^{\theta b_{j}}Y_{[2^{\nu_{1}-1}],\ldots, [2^{\nu_{m}-1}]}^{\theta}(f)_{p, \tau}\Biggr)^{\frac{1}{\theta}}   
    \end{equation}
    is valid.
  \end{theorem}
   \proof By the monotonicity property of the mixed modulus of smoothness and the relation
   \begin{equation}\label{eq2 2}
  \int_{\frac{1}{2^{\nu}}}^{\frac{1}{2^{\nu - 1}}}\frac{(1-\log t)^{\theta b_{j}}}{t}dt \asymp \nu^{\theta b_{j}},
  \end{equation}
 for 
  $b_{j}> - \frac{1}{\theta}$,  $j=1,\ldots, m$
 we have 
  \begin{multline}\label{eq2 3}
  \int_{0}^{1}...\int_{0}^{1}\omega_{\bar k}^{\theta}(f, \bar{t})_{p, \tau}\prod_{j=1}^{m} \frac{(1-\log t_{j})^{\theta b_{j}}}{t_{j}}dt_{1}...dt_{m}
  \\
  \leq \sum\limits_{\nu_{m} =1}^{\infty}...\sum\limits_{\nu_{1} =1}^{\infty}\omega_{\bar k}^{\theta}(f, \frac{1}{2^{\nu_{1} - 1}}, \ldots , \frac{1}{2^{\nu_{m} - 1}})_{p, \tau}
\int_{\frac{1}{2^{\nu_{1}}}}^{\frac{1}{2^{\nu_{1} - 1}}}...\int_{\frac{1}{2^{\nu_{m}}}}^{\frac{1}{2^{\nu_{m} - 1}}} \prod_{j=1}^{m} \frac{(1-\log t_{j})^{\theta b_{j}}}{t_{j}}dt_{1}...dt_{m} 
\\
\leq C \sum\limits_{\nu_{m} =1}^{\infty}...\sum\limits_{\nu_{1} =1}^{\infty}\prod_{j=1}^{m}\nu_{j}^{\theta b_{j}}\omega_{\bar k}^{\theta}(f, \frac{1}{2^{\nu_{1} - 1}}, \ldots , \frac{1}{2^{\nu_{m} - 1}})_{p, \tau}.
  \end{multline}
 By the property of the smoothness modulus and the ratio (2. 2), it is not difficult to make sure that
    \begin{multline}\label{eq2 4}
  \int_{0}^{1}...\int_{0}^{1}\omega_{\bar k}^{\theta}(f, \bar{t})_{p, \tau}\prod_{j=1}^{m} \frac{(1-\log t_{j})^{\theta b_{j}}}{t_{j}}dt_{1}...dt_{m}
  \\
  \geq C \sum\limits_{\nu_{m} =1}^{\infty}...\sum\limits_{\nu_{1} =1}^{\infty}\omega_{\bar k}^{\theta}(f, \frac{1}{2^{\nu_{1} - 1}}, \ldots , \frac{1}{2^{\nu_{m} - 1}})_{p, \tau}\int_{\frac{1}{2^{\nu_{1}}}}^{\frac{1}{2^{\nu_{1} - 1}}}...\int_{\frac{1}{2^{\nu_{m}}}}^{\frac{1}{2^{\nu_{m} - 1}}} \prod_{j=1}^{m} \frac{(1-\log t_{j})^{\theta b_{j}}}{t_{j}}dt_{1}...dt_{m} 
\\
\geq C \sum\limits_{\nu_{m} =1}^{\infty}...\sum\limits_{\nu_{1} =1}^{\infty}\prod_{j=1}^{m}\nu_{j}^{\theta b_{j}}\omega_{\bar k}^{\theta}(f, \frac{1}{2^{\nu_{1} - 1}}, \ldots , \frac{1}{2^{\nu_{m} - 1}})_{p, \tau},
  \end{multline}
  for  
  $b_{j}> - \frac{1}{\theta}$,  $j=1,\ldots, m$.
Now using Lemma 4, from inequality (2. 4), we get
  \begin{multline}\label{eq2 5} 
   \sum\limits_{\nu_{m} =1}^{\infty}...\sum\limits_{\nu_{1} =1}^{\infty}\prod_{j=1}^{m}(\nu_{j} + 1)^{\theta b_{j}}Y_{[2^{\nu_{1}-1}],\ldots, [2^{\nu_{m}-1}]}^{\theta}(f)_{p, \tau}
\\
   \leq C  \int_{0}^{1}...\int_{0}^{1}\omega_{\bar k}^{\theta}(f, \bar{t})_{p, \tau}\prod_{j=1}^{m} \frac{(1-\log t_{j})^{\theta b_{j}}}{t_{j}}dt_{1}...dt_{m},
   \end{multline}
 for  $b_{j}> - \frac{1}{\theta}$,  $j=1,\ldots, m$.
   
According to the inverse theorem of the `angle` approximation theory (Lemma 5), it follows from inequality (2.3) that
  \begin{multline}\label{eq2 6} 
   \int_{0}^{1}...\int_{0}^{1}\omega_{\bar k}^{\theta}(f, \bar{t})_{p, \tau}\prod_{j=1}^{m} \frac{(1-\log t_{j})^{\theta b_{j}}}{t_{j}}dt_{1}...dt_{m}
   \\
   \leq C\sum\limits_{\nu_{m} =1}^{\infty}...\sum\limits_{\nu_{1} =1}^{\infty}\prod_{j=1}^{m}\nu_{j}^{\theta b_{j}} \prod_{j=1}^{m}(2^{\nu_{j}} + 1)^{-k_{j}\theta}\Bigl(\sum\limits_{l_{m} =0}^{2^{\nu_{m}}}...\sum\limits_{\nu_{1} =0}^{2^{\nu_{m}}}\prod_{j=1}^{m}(l_{j} + 1)^{k_{j} - 1}Y_{l_{1},\ldots, l_{m}}(f)_{p, \tau} \Bigr)^{\theta},
   \end{multline}  
 for 
  $b_{j}> - \frac{1}{\theta}$,  $j=1,\ldots, m$.
Due to the monotonous decreasing of the best approximation <<angle>>, we have
\begin{multline}\label{eq2 7}
\sum\limits_{l_{m} =0}^{2^{\nu_{m}}}...\sum\limits_{\nu_{1} =0}^{2^{\nu_{m}}}\prod_{j=1}^{m}(l_{j} + 1)^{k_{j} - 1}Y_{l_{1},\ldots, l_{m}}(f)_{p, \tau}
\\
 = \sum\limits_{\mu_{m} =0}^{\nu_{m}}...\sum\limits_{\mu_{1} =0}^{\nu_{1}} \sum\limits_{l_{m} =[2^{\mu_{m} - 1}]}^{2^{\mu_{m}}}...\sum\limits_{\nu_{1} =[2^{\mu_{1} - 1}]}^{2^{\nu_{1}}}\prod_{j=1}^{m}(l_{j} + 1)^{k_{j} - 1}Y_{l_{1},\ldots, l_{m}}(f)_{p, \tau}
 \\
 \leq C \sum\limits_{\mu_{m} =0}^{\nu_{m}}...\sum\limits_{\mu_{1} =0}^{\nu_{1}}\prod_{j=1}^{m}2^{\mu_{j}k_{j}}Y_{[2^{\mu_{1} - 1}],\ldots, [2^{\mu_{m} - 1}]}(f)_{p, \tau}.
\end{multline}  
Now it follows from inequalities (2.6) and (2.7) that
\begin{multline*}
\int_{0}^{1}...\int_{0}^{1}\omega_{\bar k}^{\theta}(f, \bar{t})_{p, \tau}\prod_{j=1}^{m} \frac{(1-\log t_{j})^{\theta b_{j}}}{t_{j}}dt_{1}...dt_{m}
\\
\leq C\sum\limits_{\nu_{m} =0}^{\infty}...\sum\limits_{\nu_{1} =0}^{\infty}\prod_{j=1}^{m}\nu_{j}^{\theta b_{j}} \prod_{j=1}^{m}(2^{\nu_{j}} + 1)^{-k_{j}\theta} 
\Bigl(\sum\limits_{\mu_{m} =0}^{\nu_{m}}...\sum\limits_{\mu_{1} =0}^{\nu_{1}}\prod_{j=1}^{m}2^{\mu_{j}k_{j}}Y_{[2^{\mu_{1} - 1}],\ldots, [2^{\mu_{m} - 1}]}(f)_{p, \tau} , \Bigr)^{\theta}
\end{multline*}
 for   $b_{j}> - \frac{1}{\theta}$,  $j=1,\ldots, m$.
Since  
  \begin{equation*}
  \sum\limits_{\nu_{j} =\mu_{j}}^{\infty}(\nu_{j}+)^{\theta b_{j}}(2^{\nu_{j}} + 1)^{-k_{j}\theta}\leq C(\mu_{j}+1)^{\theta b_{j}}2^{-\mu_{j}k_{j}\theta}, 
  \end{equation*} 
for  
 $j=1,\ldots, m$, then using the generalized Hardy inequality (see \cite{7}, \cite{26}), from here, we get
  \begin{multline}\label{eq2 8}   
  \int_{0}^{1}...\int_{0}^{1}\omega_{\bar k}^{\theta}(f, \bar{t})_{p, \tau}\prod_{j=1}^{m} \frac{(1-\log t_{j})^{\theta b_{j}}}{t_{j}}dt_{1}...dt_{m}
\\  
  \ll \sum\limits_{\nu_{m} =0}^{\infty}...\sum\limits_{\nu_{1} =0}^{\infty}\prod_{j=1}^{m}(\nu_{j} + 1)^{\theta b_{j}}Y_{[2^{\mu_{1} - 1}],\ldots, [2^{\mu_{m} - 1}]}^{\theta}(f)_{p, \tau},
  \end{multline} 
for $b_{j}> - \frac{1}{\theta}$,  $j=1,\ldots, m$.  
    \hfill $\Box$
  
  \begin{theorem}\label{th2}  
Let $1 <p <+ \infty$, $1< \tau < \infty$, $0< \theta \leq \infty$ and the numbers $b_{j}> - \frac{1}{\theta}$, for $j=1,\ldots, m$. Then for the functions $f\in S_{p, \theta}^{\bar r}\mathbf{B}$ the relation 
 \begin{equation*} 
 \|f\|_{S_{p, \tau, \theta}^{0, \bar b}\mathbf{B}}
\asymp \|f\|_{p, \tau} + \Biggl(\sum\limits_{\nu_{m} =1}^{\infty}...\sum\limits_{\nu_{1} =1}^{\infty}\prod_{j=1}^{m}(\nu_{j} + 1)^{\theta b_{j}}\Biggl\|\Biggl(\sum\limits_{s_{m}=\nu_{m}}^{\infty}...\sum\limits_{s_{1}=\nu_{1}}^{\infty}|\delta_{\overline{s}}(f)|^{2}\Biggr)^{\frac{1}{2}}\Biggr\|_{p, \tau}^{\theta}\Biggr)^{\frac{1}{\theta}}    
    \end{equation*}
    is valid.
\end{theorem}
 \proof
According to the Littlewood--Paley theorem in the Lorentz space (see for example \cite{24}) for the function $f\in L_{p, \tau}(\mathbb{T}^{m})$ there is an inequality
\begin{multline}\label{eq2 9}
Y_{[2^{\nu_{1} - 1}],\ldots, [2^{\nu_{m} - 1}]}(f)_{p, \tau} \leq \|f - U_{[2^{\nu_{1} - 1}],\ldots, [2^{\nu_{m} - 1}]}(f)\|_{p, \tau}
\\
=\Bigl\|\sum\limits_{s_{m}=\nu_{m}}^{\infty}...\sum\limits_{s_{1}=\nu_{1}}^{\infty}\sum\limits_{\overline{n} \in \rho \left( \overline{s} \right)
}a_{\overline{n} }(f) e^{i\langle\overline{n} ,\overline{x}\rangle}\Bigr\|_{p, \tau} \asymp \Biggl\|\Biggl(\sum\limits_{s_{m}=\nu_{m}}^{\infty}...\sum\limits_{s_{1}=\nu_{1}}^{\infty}|\delta_{\overline{s}}(f)|^{2}\Biggr)^{\frac{1}{2}}\Biggr\|_{p, \tau},
\end{multline}
for  $1< p < \infty$, $1< \tau < \infty$.
Now  from inequalities (2.8) and (2.9)  it follows that
 \begin{multline}\label{eq2 10} 
   \int_{0}^{1}...\int_{0}^{1}\omega_{\bar k}^{\theta}(f, \bar{t})_{p, \tau}\prod_{j=1}^{m} \frac{(1-\log t_{j})^{\theta b_{j}}}{t_{j}}dt_{1}...dt_{m}
   \\
   \leq C\sum\limits_{\nu_{m} =0}^{\infty}...\sum\limits_{\nu_{1} =0}^{\infty}\prod_{j=1}^{m}(\nu_{j} + 1)^{\theta b_{j}}\Biggl\|\Biggl(\sum\limits_{s_{m}=\nu_{m}}^{\infty}...\sum\limits_{s_{1}=\nu_{1}}^{\infty}|\delta_{\overline{s}}(f)|^{2}\Biggr)^{\frac{1}{2}}\Biggr\|_{p, \tau}^{\theta} 
\end{multline}
for  $b_{j}> - \frac{1}{\theta}$,  $j=1,\ldots, m$, $1< p < \infty$.
 
 Let's prove the inverse inequality. It is known that
  \begin{equation*}
 \|f - U_{[2^{\nu_{1} - 1}],\ldots, [2^{\nu_{m} - 1}]}(f)\|_{p, \tau} \leq CY_{[2^{\nu_{1} - 1}],\ldots, [2^{\nu_{m} - 1}]}(f)_{p, \tau},
 \end{equation*}
for the function $f\in L_{p, \tau}(\mathbb{T}^{m})$, for $1<p<\infty$. 
 Therefore,  from inequality (2. 5) and relation (2. 9) it follows that
  \begin{multline}\label{eq2 11}
\sum\limits_{\nu_{m} =1}^{\infty}...\sum\limits_{\nu_{1} =1}^{\infty}\prod_{j=1}^{m}\nu_{j}^{\theta b_{j}}\Biggl\|\Biggl(\sum\limits_{s_{m}=\nu_{m}}^{\infty}...\sum\limits_{s_{1}=\nu_{1}}^{\infty}|\delta_{\overline{s}}(f)|^{2}\Biggr)^{\frac{1}{2}}\Biggr\|_{p, \tau}^{\theta} 
\\
\leq C \int_{0}^{1}...\int_{0}^{1}\omega_{\bar k}^{\theta}(f, \bar{t})_{p, \tau}\prod_{j=1}^{m} \frac{(1-\log t_{j})^{\theta b_{j}}}{t_{j}}dt_{1}...dt_{m},
\end{multline}
 for   $b_{j}> - \frac{1}{\theta}$,  $j=1,\ldots, m$, $1< p < \infty$.
Further considering that
\begin{equation*} 
\Biggl\|\Biggl(\sum\limits_{s_{1}=0}^{\infty}...\sum\limits_{s_{j}=0}^{\infty}\sum\limits_{s_{j+1}=\nu_{j+1}}^{\infty}\ldots \sum\limits_{s_{m}=\nu_{m}}^{\infty}|\delta_{\overline{s}}(f)|^{2}\Biggr)^{\frac{1}{2}}\Biggr\|_{p, \tau}
\leq \Biggl\|\Biggl(\sum\limits_{s_{1}=0}^{\infty}\ldots \sum\limits_{s_{m}=0}^{\infty}|\delta_{\overline{s}}(f)|^{2}\Biggr)^{\frac{1}{2}}\Biggr\|_{p}\asymp \|f\|_{p, \tau},
\end{equation*} 
 for the function $f\in L_{p, \tau}(\mathbb{T}^{m})$, for $1<p <\infty$, $1< \tau< \infty$, the inequality (2. 10) can be rewritten in the following form
\begin{multline}\label{eq2 12} 
   \int_{0}^{1}...\int_{0}^{1}\omega_{\bar k}^{\theta}(f, \bar{t})_{p, \tau}\prod_{j=1}^{m} \frac{(1-\log t_{j})^{\theta b_{j}}}{t_{j}}dt_{1}...dt_{m}
   \\
   \leq C\Biggl\{\|f\|_{p, \tau}^{\theta} +\sum\limits_{\nu_{m} =1}^{\infty}...\sum\limits_{\nu_{1} =1}^{\infty}\prod_{j=1}^{m}(\nu_{j} + 1)^{\theta b_{j}}\Biggl\|\Biggl(\sum\limits_{s_{m}=\nu_{m}}^{\infty}...\sum\limits_{s_{1}=\nu_{1}}^{\infty}|\delta_{\overline{s}}(f)|^{2}\Biggr)^{\frac{1}{2}}\Biggr\|_{p, \tau}^{\theta}\Biggr\} 
\end{multline}
for $b_{j}> - \frac{1}{\theta}$,  $j=1,\ldots, m$, $1< p, \tau < \infty$.

Now, from inequalities (2. 11) and (2. 12), by definition of the class, we obtain the statement of Theorem 2.
   \hfill $\Box$

\begin{theorem}\label{th3}
Let $1 <p, \tau <+ \infty$, $0< \theta \leq \infty$ and the numbers $b_{j}> - \frac{1}{\theta}$, for $j=1,\ldots, m$. Then for the function $f\in S_{p, \tau, \theta}^{0, \bar b}\mathbf{B}$, holds the following relation: 
 \begin{multline*} 
C_{1} \Biggl\{\|f\|_{p, \tau} + \Biggl(\sum\limits_{l_{m} =1}^{\infty}...\sum\limits_{l_{1} =1}^{\infty}\prod_{j=1}^{m}2^{l_{j}\theta (b_{j} + \frac{1}{\theta})}\Biggl\|\sum\limits_{s_{m}=[2^{l_{m}-1}]+1}^{2^{l_{m}}}...\sum\limits_{s_{1}=[2^{l_{1}-1}]+1}^{2^{l_{1}}}\delta_{\overline{s}}(f)\Biggr\|_{p, \tau}^{\theta}\Biggr)^{\frac{1}{\theta}}\Biggr\}
\leq \|f\|_{S_{p, \tau, \theta}^{0, \bar b}\mathbf{B}}
\\
\leq C_{2}\Biggl\{ \|f\|_{p, \tau} + \Biggl(\sum\limits_{l_{m} =0}^{\infty}...\sum\limits_{l_{1} =0}^{\infty}\prod_{j=1}^{m}2^{l_{j}\theta (b_{j}+ \frac{1}{\theta})}\Biggl\|\sum\limits_{s_{m}=[2^{l_{m}-1}]+1}^{2^{l_{m}}}...\sum\limits_{s_{1}=[2^{l_{1}-1}]+1}^{2^{l_{1}}}\delta_{\overline{s}}(f)\Biggr\|_{p, \tau}^{\theta}\Biggr)^{\frac{1}{\theta}}\Biggr\}.    
    \end{multline*}
  \end{theorem}

\proof 
Since
\begin{equation*}
\sum\limits_{\nu=[2^{l1}]+1}^{2^{l}}\nu^{b_{j}\theta}\asymp 2^{l\theta(b_{j}+\frac{1}{\theta})}, \, \, j=1, \ldots, m
\end{equation*}
and the sequence
\begin{equation*}
\sigma_{\nu_{1},...,\nu_{m}}=\Biggl\|\Biggl(\sum\limits_{s_{m}=\nu_{m}}^{\infty}...\sum\limits_{s_{1}=\nu_{1}}^{\infty}|\delta_{\overline{s}}(f)|^{2}\Biggr)^{\frac{1}{2}}\Biggr\|_{p, \tau}
\end{equation*}
monotonically decreases for each index $\nu_{j}$, for $j=1, \ldots, m$, then
\begin{multline}\label{eq2 13}
\sum\limits_{\nu_{m} =1}^{\infty}...\sum\limits_{\nu_{1} =1}^{\infty}\prod_{j=1}^{m}(\nu_{j} + 1)^{\theta b_{j}}\Biggl\|\Biggl(\sum\limits_{s_{m}=\nu_{m}}^{\infty}...\sum\limits_{s_{1}=\nu_{1}}^{\infty}|\delta_{\overline{s}}(f)|^{2}\Biggr)^{\frac{1}{2}}\Biggr\|_{p, \tau}^{\theta}= \sum\limits_{\nu_{m} =1}^{\infty}...\sum\limits_{\nu_{1} =1}^{\infty}\prod_{j=1}^{m}(\nu_{j} + 1)^{\theta b_{j}}\sigma_{\nu_{1},...,\nu_{m}}^{\theta}
\\
\leq C\sum\limits_{l_{m} =0}^{\infty}...\sum\limits_{l_{1} =0}^{\infty}\prod_{j=1}^{m}2^{l_{j}\theta(b_{j} + \frac{1}{\theta})} \sigma_{[2^{l_{1}-1}]+1,...,[2^{l_{m}-1}]+1}^{\theta}.
\end{multline}
By the Littlewood--Paley theorem \cite{24} and the triangle inequality for the norm we will have
\begin{multline*} 
\sigma_{[2^{l_{1}-1}]+1,...,[2^{l_{m}-1}]+1} \leq C\Biggl\|\sum\limits_{s_{m}=[2^{l_{m}-1}]+1}^{\infty}...\sum\limits_{s_{1}=[2^{l_{1}-1}]+1}^{\infty}\delta_{\overline{s}}(f)\Biggr\|_{p, \tau} 
\\
\leq C \sum\limits_{\mu_{m}=l_{m}}^{\infty}...\sum\limits_{\mu_{1}=l_{1}}^{\infty}\Biggl\|\sum\limits_{s_{m}=[2^{\mu_{m}-1}]+1}^{2^{\mu_{m}}}...\sum\limits_{s_{1}=[2^{\mu_{1}-1}]+1}^{2^{\mu_{1}}}\delta_{\overline{s}}(f)\Biggr\|_{p, \tau}.
\end{multline*}
Therefore , from (2. 13) it follows that
 \begin{multline}\label{eq2 14}
\sum\limits_{\nu_{m} =1}^{\infty}...\sum\limits_{\nu_{1} =1}^{\infty}\prod_{j=1}^{m}(\nu_{j} + 1)^{\theta b_{j}}\Biggl\|\Biggl(\sum\limits_{s_{m}=\nu_{m}}^{\infty}...\sum\limits_{s_{1}=\nu_{1}}^{\infty}|\delta_{\overline{s}}(f)|^{2}\Biggr)^{\frac{1}{2}}\Biggr\|_{p, \tau}^{\theta} 
\\
\leq C \sum\limits_{l_{m} =0}^{\infty}...\sum\limits_{l_{1} =0}^{\infty}\prod_{j=1}^{m}2^{l_{j}\theta(b_{j} + \frac{1}{\theta})} \Biggl(\sum\limits_{\mu_{m}=l_{m}}^{\infty}...\sum\limits_{\mu_{1}=l_{1}}^{\infty}\Biggl\|\sum\limits_{s_{m}=[2^{\mu_{m}-1}]+1}^{2^{\mu_{m}}}...\sum\limits_{s_{1}=[2^{\mu_{1}-1}]+1}^{2^{\mu_{1}}}\delta_{\overline{s}}(f)\Biggr\|_{p, \tau}\Biggr)^{\theta}.
\end{multline} 
Since $b_{j} + \frac{1}{\theta}> 0$ for $j=1, \ldots, m$, then
$\sum\limits_{l_{j} =0}^{\mu_{j}}2^{l_{j}\theta(b_{j} + \frac{1}{\theta})}\leq C2^{\mu_{j}\theta(b_{j} + \frac{1}{\theta})}, \, \,  j=1, \ldots, m.
$
Therefore, using the generalized Hardy inequality \cite{26}, from the formula (2. 14) we obtain
\begin{multline}\label{eq2 15}
\sum\limits_{\nu_{m} =1}^{\infty}...\sum\limits_{\nu_{1} =1}^{\infty}\prod_{j=1}^{m}(\nu_{j} + 1)^{\theta b_{j}}\Biggl\|\Biggl(\sum\limits_{s_{m}=\nu_{m}}^{\infty}...\sum\limits_{s_{1}=\nu_{1}}^{\infty}|\delta_{\overline{s}}(f)|^{2}\Biggr)^{\frac{1}{2}}\Biggr\|_{p, \tau}^{\theta}
\\
\leq C\sum\limits_{l_{m} =0}^{\infty}...\sum\limits_{l_{1} =0}^{\infty}\prod_{j=1}^{m}2^{l_{j}\theta(b_{j} + \frac{1}{\theta})}\Biggl\|\sum\limits_{s_{m}=[2^{l_{m}-1}]+1}^{2^{l_{m}}}...\sum\limits_{s_{1}=[2^{l_{1}-1}]+1}^{2^{l_{1}}}\delta_{\overline{s}}(f)\Biggr\|_{p, \tau}^{\theta}.
\end{multline} 
Now from inequality (2. 15) and Theorem 2 follows the second inequality in Theorem 3.
 
 To prove the first inequality in Theorem 3, according to
of the monotone decreasing sequence $\{\sigma_{\nu_{1},...,\nu_{m}}\}$ and the Littlewood--Paley theorem \cite{24} we have  
\begin{multline}\label{eq2 16}
\sum\limits_{\nu_{m} =1}^{\infty}...\sum\limits_{\nu_{1} =1}^{\infty}\prod_{j=1}^{m}(\nu_{j} + 1)^{\theta b_{j}}\sigma_{\nu_{1},...,\nu_{m}}^{\theta}\geq C\sum\limits_{l_{m} =1}^{\infty}...\sum\limits_{l_{1} =1}^{\infty}\prod_{j=1}^{m}2^{l_{j}\theta(b_{j} + \frac{1}{\theta})} \sigma_{[2^{l_{1}-1}]+1,...,[2^{l_{m}-1}]+1}^{\theta} 
\\
\geq C\sum\limits_{l_{m} =1}^{\infty}...\sum\limits_{l_{1} =1}^{\infty}\prod_{j=1}^{m}2^{l_{j}\theta(b_{j} + \frac{1}{\theta})}\Biggl\|\Biggl(\sum\limits_{s_{m}=2^{l_{m}-1}+1}^{2^{l_{m}}}...\sum\limits_{s_{1}=2^{l_{1}-1}+1}^{2^{l_{1}}}|\delta_{\overline{s}}(f)|^{2}\Biggr)^{\frac{1}{2}}\Biggr\|_{p, \tau}^{\theta} 
\\
\geq C\sum\limits_{l_{m} =1}^{\infty}...\sum\limits_{l_{1} =1}^{\infty}\prod_{j=1}^{m}2^{l_{j}\theta(b_{j} + \frac{1}{\theta})}\Biggl\|\sum\limits_{s_{m}=2^{l_{m}-1}+1}^{2^{l_{m}}}...\sum\limits_{s_{1}=2^{l_{1}-1}+1}^{2^{l_{1}}}\delta_{\overline{s}}(f)\Biggr\|_{p, \tau}^{\theta}.
\end{multline}
Now from inequality (2. 16) and Theorem 2 follows the first inequality in Theorem 3.
      \hfill $\Box$

In the following theorems, we establish the conditions for embedding classes $S_{p,\theta}^{0, \overline{b}}\mathbf{B}$ and $S_{p, \theta}^{0, \overline{b}}B$.
\begin{theorem}\label{th4}   
  Let $1 <p, \to <+ \infty$, $0< \theta\leq\infty$ and the numbers $b_{j}> - \frac{1}{\theta}$, for $j=1,\ldots, m$. Then
    $S_{p, \tau, \theta}^{0, \overline{v}}B \subset S_{p, \tau, \theta}^{0, \overline{b}}\mathbf{B} \subset S_{p, \tau, \theta}^{0, \overline{u}}B$, where $\overline{v}=(v_{1}, \ldots ,v_{m})$, $\overline{u}=(u_{1}, \ldots ,u_{m})$, $v_{j}=b_{j}+\frac{1}{\min\{\beta, \theta\}}$, $u_{j}=b_{j}+\frac{1}{\max\{\gamma, \theta\}}$, for $j=1,\ldots, m$, where 
    \begin{equation*}
\beta =\left\{\begin{array}{rl} \tau, & \mbox{if} \, \,  1<\tau \leq 2 \, \, \mbox{and} \, \, 1<p< \infty \\
 2, & \mbox{if} \, \, 2< \tau < \infty \, \,  \mbox{and} \, \, 2<p< \infty ;
\end{array}
 \right.    
   \end{equation*} 
  \begin{equation*}
\gamma =\left\{\begin{array}{rl} 2, & \mbox{if} \, \,  1<\tau \leq 2 \, \, \mbox{and} \, \, 1<p< \infty \\
 \tau, & \mbox{if} \, \, 2< \tau < \infty \, \,  \mbox{and} \, \, 2<p< \infty .
\end{array}
 \right.    
   \end{equation*}          
\end{theorem}   
 
 \proof
For the function $f\in L_{p, \tau}(\mathbb{T}^{m})$, for $1<p< \infty$ it is known that (see for example \cite{22}) 
 \begin{equation}\label{eq2 17} 
 \Biggl\|\Biggl(\sum\limits_{s_{m}=\nu_{m}}^{\infty}...\sum\limits_{s_{1}=\nu_{1}}^{\infty}|\delta_{\overline{s}}(f)|^{2}\Biggr)^{\frac{1}{2}}\Biggr\|_{p, \tau}\leq C \Biggl(\sum\limits_{s_{m}=\nu_{m}}^{\infty}...\sum\limits_{s_{1}=\nu_{1}}^{\infty}\|\delta_{\overline{s}}(f)\|_{p, \tau}^{\beta}\Biggr)^{\frac{1}{\beta}}.
   \end{equation}
 
If $0< \theta\leq\beta$, then according to Jensen's inequality ((\cite[ch.~3, sec.~3, p.~125]{2}) from the formula (2. 17) we get
 \begin{equation*} 
 \Biggl\|\Biggl(\sum\limits_{s_{m}=\nu_{m}}^{\infty}...\sum\limits_{s_{1}=\nu_{1}}^{\infty}|\delta_{\overline{s}}(f)|^{2}\Biggr)^{\frac{1}{2}}\Biggr\|_{p, \tau}^{\theta}\leq C \sum\limits_{s_{m}=\nu_{m}}^{\infty}...\sum\limits_{s_{1}=\nu_{1}}^{\infty}\|\delta_{\overline{s}}(f)\|_{p, \tau}^{^{\theta}}.
   \end{equation*}
Therefore , by Theorem 2 we have
\begin{equation}\label{eq2 18}   
   \|f\|_{S_{p, \tau, \theta}^{0, \overline{b}}\mathbf{B}} \leq C\Bigl\{\|f\|_{p, \tau} + \Biggl(\sum\limits_{\nu_{m} =1}^{\infty}...\sum\limits_{\nu_{1} =1}^{\infty}\prod_{j=1}^{m}(\nu_{j} + 1)^{\theta b_{j}}\sum\limits_{s_{m}=\nu_{m}}^{\infty}...\sum\limits_{s_{1}=\nu_{1}}^{\infty}\|\delta_{\overline{s}}(f)\|_{p, \tau}^{^{\theta}}\Biggr)^{\frac{1}{\theta}} \Bigr\}. 
   \end{equation}
  Further changing the summation order and considering that $b_{j}> is \frac{1}{\theta}$, for $j=1,\ldots, m$ we get
   \begin{multline}\label{eq2 19}
   \sum\limits_{\nu_{m} =1}^{\infty}...\sum\limits_{\nu_{1} =1}^{\infty}\prod_{j=1}^{m}(\nu_{j} + 1)^{\theta b_{j}}\sum\limits_{s_{m}=\nu_{m}}^{\infty}...\sum\limits_{s_{1}=\nu_{1}}^{\infty}\|\delta_{\overline{s}}(f)\|_{p, \tau}^{^{\theta}} 
   \\
   = \sum\limits_{s_{m} =1}^{\infty}...\sum\limits_{s_{1} =1}^{\infty}\|\delta_{\overline{s}}(f)\|_{p, \tau}^{^{\theta}}\sum\limits_{\nu_{m} =1}^{s_{m}}...\sum\limits_{\nu_{1} =1}^{s_{1}}\prod_{j=1}^{m}(\nu_{j} + 1)^{\theta b_{j}}\leq C 
\sum\limits_{s_{m} =1}^{\infty}...\sum\limits_{s_{1} =1}^{\infty}\prod_{j=1}^{m}s_{j} ^{\theta b_{j}+1}\|\delta_{\overline{s}}(f)\|_{p, \tau}^{^{\theta}}.  
   \end{multline} 
Now  from inequalities (2. 18) and (2. 19) it follows that
   $S_{p, \tau, \theta}^{0, \overline{v}}B \subset S_{p, \tau, \theta}^{0, \overline{b}}\mathbf{B}$,  in the case $0< \theta \leq \beta$.

Let 
 $1< \beta < \theta \leq \infty$. Then applying the H\"{o}lder inequality for $\eta=\frac{\beta}{\beta}$ and $\frac{1}{\eta}+\frac{1}{\beta^{'}}=1$ from (2. 17) we get
   \begin{multline}\label{eq2 20} 
 \Biggl\|\Biggl(\sum\limits_{s_{m}=\nu_{m}}^{\infty}...\sum\limits_{s_{1}=\nu_{1}}^{\infty}|\delta_{\overline{s}}(f)|^{2}\Biggr)^{\frac{1}{2}}\Biggr\|_{p, \tau}\leq C \Biggl(\sum\limits_{s_{m}=\nu_{m}}^{\infty}...\sum\limits_{s_{1}=\nu_{1}}^{\infty}\prod_{j=1}^{m}s_{j} ^{\theta b_{j}}\|\delta_{\overline{s}}(f)\|_{p, \tau}^{\theta}\Biggr)^{\frac{1}{\theta}}
\\ 
 \times
 \Biggl(\sum\limits_{s_{m}=\nu_{m}}^{\infty}...\sum\limits_{s_{1}=\nu_{1}}^{\infty}\prod_{j=1}^{m}s_{j}^{-b_{j}\beta\eta^{'}} \Biggr)^{\frac{1}{\beta\eta^{'}}}\leq C \Biggl(\sum\limits_{s_{m}=\nu_{m}}^{\infty}...\sum\limits_{s_{1}=\nu_{1}}^{\infty}\prod_{j=1}^{m}s_{j} ^{\theta b_{j}}\|\delta_{\overline{s}}(f)\|_{p, \tau}^{\theta}\Biggr)^{\frac{1}{\theta}}\prod_{j=1}^{m}\nu_{j}^{\frac{1}{\beta\eta^{'}}-b_{j}}
   \end{multline}
   for 
    $b_{j}> \frac{1}{\beta}-\frac{1}{\theta}$,  $j=1,\ldots, m$.
Now  from (2. 20) and Theorem 2 it follows that
    \begin{equation}\label{eq2 21}
    \|f\|_{S_{p, \tau, \theta}^{0, \overline{b}}\mathbf{B}} 
    \leq C\Bigl\{\|f\|_{p, \tau} + \Biggl(\sum\limits_{\nu_{m} =1}^{\infty}...\sum\limits_{\nu_{1} =1}^{\infty}\prod_{j=1}^{m}(\nu_{j} + 1)^{\frac{\theta}{\beta} -1}\sum\limits_{s_{m}=\nu_{m}}^{\infty}...\sum\limits_{s_{1}=\nu_{1}}^{\infty}\prod_{j=1}^{m}s_{j} ^{\theta b_{j}}\|\delta_{\overline{s}}(f)\|_{p, \tau}^{\theta}\Biggr)^{\frac{1}{\theta}} \Bigr\},
   \end{equation} 
   $1< \beta < \theta \leq \infty$, $b_{j}> \frac{1}{\beta}-\frac{1}{\theta}$, for $j=1,\ldots, m$.
   
   Next, changing the order of summation and taking into account, that
$\sum\limits_{\nu=1}^{s}\nu^{\frac{\theta}{\beta} -1}\leq C s^{\frac{\theta}{\beta}},
$ from inequality (2. 21) we get
 \begin{multline}\label{eq2 22}
    \|f\|_{S_{p, \tau, \theta}^{0, \bar b}\mathbf{B}} \leq C\Biggl\{\|f\|_{p, \tau} + \Biggl(\sum\limits_{s_{m} =1}^{\infty}...\sum\limits_{s_{1} =1}^{\infty}\prod_{j=1}^{m}s_{j} ^{\theta b_{j}}\|\delta_{\overline{s}}(f)\|_{p, \tau}^{^{\theta}}\sum\limits_{\nu_{m}=1}^{s_{m}}...\sum\limits_{\nu_{1}=1}^{s_{1}}\prod_{j=1}^{m}(\nu_{j} + 1)^{\frac{\theta}{\beta} -1}\Biggr)^{\frac{1}{\theta}} \Biggr\}
    \\
    \leq C\Biggl\{\|f\|_{p, \tau} + \Biggl(\sum\limits_{s_{m} =1}^{\infty}...\sum\limits_{s_{1} =1}^{\infty}\prod_{j=1}^{m}s_{j} ^{\theta (b_{j}+\frac{1}{\beta})}\|\delta_{\overline{s}}(f)\|_{p, \tau}^{^{\theta}}\Biggr)^{\frac{1}{\theta}}\Biggr\},
   \end{multline}
 in the case  
     $1< \beta < \theta \leq \infty$, $b_{j}> \frac{1}{\beta}-\frac{1}{\theta}$, for $j=1,\ldots, m$. Therefore 
      $S_{p, \tau, \theta}^{0, \overline{v}}B \subset S_{p, \tau, \theta}^{0, \overline{b}}\mathbf{B}$,  in the case 
          $1< \beta < \theta \leq \infty$, $b_{j}> \frac{1}{\beta}-\frac{1}{\theta}$, for $j=1,\ldots, m$.

We prove the inclusion of $S_{p, \tau, \theta}^{0, \overline{b}}\mathbf{B}\subset_{p, \tau, \theta}^{0, \overline{u}}B$.
 For the function $f\in L_{p, \tau}(\mathbb{T}^{m})$, for $1<p< \infty$ it is known that (see for example \cite{24})
\begin{equation}\label{eq2 23} 
 \Biggl\|\Biggl(\sum\limits_{s_{m}=\nu_{m}}^{\infty}...\sum\limits_{s_{1}=\nu_{1}}^{\infty}|\delta_{\overline{s}}(f)|^{2}\Biggr)^{\frac{1}{2}}\Biggr\|_{p}\geq C \Biggl(\sum\limits_{s_{m}=\nu_{m}}^{\infty}...\sum\limits_{s_{1}=\nu_{1}}^{\infty}\|\delta_{\overline{s}}(f)\|_{p}^{\gamma}\Biggr)^{\frac{1}{\gamma}},
   \end{equation}
where   $\gamma=2$, if $1< \tau \leq 2$ and $1< p<\infty$ or $2< \tau < \infty$ and $2< p<\infty$.
\begin{multline}\label{eq2 24}
\Biggl(\sum\limits_{\nu_{m} =1}^{\infty}...\sum\limits_{\nu_{1} =1}^{\infty}\prod_{j=1}^{m}(\nu_{j} + 1)^{\theta b_{j}}\Biggl\|\Biggl(\sum\limits_{s_{m}=\nu_{m}}^{\infty}...\sum\limits_{s_{1}=\nu_{1}}^{\infty}|\delta_{\overline{s}}(f)|^{2}\Biggr)^{\frac{1}{2}}\Biggr\|_{p, \tau}^{\theta}\Biggr)^{\frac{1}{\theta}}
\\
\geq C\Biggl(\sum\limits_{\nu_{m} =1}^{\infty}...\sum\limits_{\nu_{1} =1}^{\infty}\prod_{j=1}^{m}(\nu_{j} + 1)^{\theta b_{j}}\Biggl(\sum\limits_{s_{m}=\nu_{m}}^{\infty}...\sum\limits_{s_{1}=\nu_{1}}^{\infty}\|\delta_{\overline{s}}(f)\|_{p, \tau}^{\gamma}\Biggr)^{\frac{\theta}{\gamma}}\Biggr)^{\frac{1}{\theta}}.
\end{multline}

If $1< \gamma\l\theta<\infty$, according to Jensen's inequality (\cite[ch.~3, sec.~3, p.~125]{2} we have
\begin{multline}\label{eq2 25}
\sum\limits_{\nu_{m} =1}^{\infty}...\sum\limits_{\nu_{1} =1}^{\infty}\prod_{j=1}^{m}(\nu_{j} + 1)^{\theta b_{j}}\Biggl(\sum\limits_{s_{m}=\nu_{m}}^{\infty}...\sum\limits_{s_{1}=\nu_{1}}^{\infty}\|\delta_{\overline{s}}(f)\|_{p, \tau}^{\gamma}\Biggr)^{\frac{\theta}{\gamma}}
\\
\geq C \sum\limits_{\nu_{m} =1}^{\infty}...\sum\limits_{\nu_{1} =1}^{\infty}\prod_{j=1}^{m}\nu_{j}^{\theta b_{j}}\sum\limits_{s_{m}=\nu_{m}}^{\infty}...\sum\limits_{s_{1}=\nu_{1}}^{\infty}\|\delta_{\overline{s}}(f)\|_{p, \tau}^{\theta} .
\end{multline}
Now changing the summation order and considering that
\begin{equation}\label{eq2 26}
\sum\limits_{\nu=1}^{s}\nu^{\theta\beta}\asymp s^{\theta\beta + 1},
\end{equation}
for $\theta\beta + 1 > 0$, we get
\begin{multline}\label{eq2 27}
\sum\limits_{\nu_{m} =1}^{\infty}...\sum\limits_{\nu_{1} =1}^{\infty}\prod_{j=1}^{m}\nu_{j}^{\theta b_{j}}\sum\limits_{s_{m}=\nu_{m}}^{\infty}...\sum\limits_{s_{1}=\nu_{1}}^{\infty}\|\delta_{\overline{s}}(f)\|_{p, \tau}^{\theta} =\sum\limits_{s_{m} =1}^{\infty}...\sum\limits_{s_{1} =1}^{\infty}\|\delta_{\overline{s}}(f)\|_{p}^{\theta}\sum\limits_{\nu_{m}=1}^{s_{m}}...\sum\limits_{\nu_{1}=1}^{s_{1}}\prod_{j=1}^{m}\nu_{j}^{\theta b_{j}}
\\
 \geq C\sum\limits_{s_{m} =1}^{\infty}...\sum\limits_{s_{1} =1}^{\infty}\prod_{j=1}^{m}s_{j}^{\theta b_{j} + 1}\|\delta_{\overline{s}}(f)\|_{p, \tau}^{\theta},  
\end{multline}
for $\theta b_{j} + 1 > 0$, $j=1,\ldots, m$.
Further, from the inequalities (2. 24), (2. 25) and (2. 27) it follows that
\begin{multline}\label{eq2 28}
\Biggl(\sum\limits_{\nu_{m} =1}^{\infty}...\sum\limits_{\nu_{1} =1}^{\infty}\prod_{j=1}^{m}(\nu_{j} + 1)^{\theta b_{j}}\Biggl\|\Biggl(\sum\limits_{s_{m}=\nu_{m}}^{\infty}...\sum\limits_{s_{1}=\nu_{1}}^{\infty}|\delta_{\overline{s}}(f)|^{2}\Biggr)^{\frac{1}{2}}\Biggr\|_{p, \tau}^{\theta}\Biggr)^{\frac{1}{\theta}}
\\
\geq C \Biggl(\sum\limits_{s_{m} =1}^{\infty}...\sum\limits_{s_{1} =1}^{\infty}\prod_{j=1}^{m}s_{j}^{\theta b_{j} + 1}\|\delta_{\overline{s}}(f)\|_{p, \tau}^{\theta}\Biggr)^{\frac{1}{\theta}}, 
\end{multline}
in the case  $1< \gamma \leq \theta< \infty$, $\theta b_{j} + 1 > 0$, $j=1,\ldots, m$.
Therefore, according to Theorem 2 of (2. 28) , we obtain that   
$S_{p, \tau, \theta}^{0, \overline{b}}\mathbf{B} \subset S_{p, \tau, \theta}^{0, \overline{u}}B$, in the case  $1< \gamma \leq \theta< \infty$, where $u_{j}=b_{j}+\frac{1}{\theta}$, for $j=1,\ldots, m$.

Let $0<\theta< \gamma$. Then $\max\{\gamma, \theta\}=\gamma$. Using the ratio (2. 26) and changing the order of summation we have
\begin{multline}\label{eq2 29}
\Biggl(\sum\limits_{s_{m} =1}^{\infty}...\sum\limits_{s_{1} =1}^{\infty}\prod_{j=1}^{m}s_{j}^{(b_{j} + \frac{1}{\gamma})\theta}\|\delta_{\overline{s}}(f)\|_{p, \tau}^{\theta}\Biggr)^{\frac{1}{\theta}} = \Biggl(\sum\limits_{s_{m} =1}^{\infty}...\sum\limits_{s_{1} =1}^{\infty}\prod_{j=1}^{m}s_{j}^{(b_{j} + \frac{1}{\gamma} + \varepsilon)\theta}\prod_{j=1}^{m}s_{j}^{-\varepsilon\theta}\|\delta_{\overline{s}}(f)\|_{p, \tau}^{\theta}\Biggr)^{\frac{1}{\theta}}  
\\
\leq C\Biggl(\sum\limits_{s_{m} =1}^{\infty}...\sum\limits_{s_{1} =1}^{\infty} \prod_{j=1}^{m}s_{j}^{-\varepsilon\theta}\|\delta_{\overline{s}}(f)\|_{p, \tau}^{\theta}\sum\limits_{\nu_{m} =1}^{s_{m}}...\sum\limits_{\nu_{1} =1}^{s_{1}}\prod_{j=1}^{m}\nu_{j}^{(b_{j} + \frac{1}{\gamma} + \varepsilon)\theta - 1} \Biggr)^{\frac{1}{\theta}}  
\\
= C \Biggl(\sum\limits_{\nu_{m} =1}^{\infty}...\sum\limits_{\nu_{1} =1}^{\infty}\prod_{j=1}^{m}\nu_{j}^{(b_{j} + \frac{1}{\gamma} + \varepsilon)\theta - 1}
\sum\limits_{s_{m}=\nu_{m}}^{\infty}...\sum\limits_{s_{1}=\nu_{1}}^{\infty}\prod_{j=1}^{m}s_{j}^{-\varepsilon\theta}\|\delta_{\overline{s}}(f)\|_{p, \tau}^{\theta}
\Biggr)^{\frac{1}{\theta}},
\end{multline}
for the number $\varepsilon > \frac{1}{\theta} - \frac{1}{\gamma}> 0$.
Next, applying the H\"{o}lder inequality for $\eta=\frac{\gamma}{\theta}$, $\frac{1}{\eta} + \frac{1}{\beta^{'}}=1$ we get

\begin{multline}\label{eq2 30}
\Biggl(\sum\limits_{s_{m}=\nu_{m}}^{\infty}...\sum\limits_{s_{1}=\nu_{1}}^{\infty}\prod_{j=1}^{m}s_{j}^{-\varepsilon\theta}\|\delta_{\overline{s}}(f)\|_{p, \tau}^{\theta}
\Biggr)^{\frac{1}{\theta}}\leq \Biggl(\sum\limits_{s_{m}=\nu_{m}}^{\infty}...\sum\limits_{s_{1}=\nu_{1}}^{\infty}\|\delta_{\overline{s}}(f)\|_{p, \tau}^{\gamma}
\Biggr)^{\frac{1}{\gamma}}
\\
\times
\Biggl(\sum\limits_{s_{m}=\nu_{m}}^{\infty}...\sum\limits_{s_{1}=\nu_{1}}^{\infty}\prod_{j=1}^{m}s_{j}^{-\varepsilon\theta\eta^{'}}\Biggr)^{\frac{1}{\theta\eta^{'}}}\leq C\prod_{j=1}^{m}\nu_{j}^{\frac{1}{\theta} - \frac{1}{\gamma}-\varepsilon}\Biggl(\sum\limits_{s_{m}=\nu_{m}}^{\infty}...\sum\limits_{s_{1}=\nu_{1}}^{\infty}\|\delta_{\overline{s}}(f)\|_{p, \tau}^{\gamma}
\Biggr)^{\frac{1}{\gamma}},
\end{multline}
for the number  $\varepsilon > \frac{1}{\theta} - \frac{1}{\gamma}> 0$.

Now,  from inequalities (2. 29) and (2. 30), it follows that
\begin{multline}\label{eq2 31}
\Biggl(\sum\limits_{s_{m} =1}^{\infty}...\sum\limits_{s_{1} =1}^{\infty}\prod_{j=1}^{m}s_{j}^{(b_{j} + \frac{1}{\gamma})\theta}\|\delta_{\overline{s}}(f)\|_{p, \tau}^{\theta}\Biggr)^{\frac{1}{\theta}} 
\\
\leq C  \Biggl(\sum\limits_{\nu_{m} =1}^{\infty}...\sum\limits_{\nu_{1} =1}^{\infty}\prod_{j=1}^{m}\nu_{j}^{(b_{j} + \frac{1}{\gamma} + \varepsilon)\theta - 1}\prod_{j=1}^{m}\nu_{j}^{(\frac{1}{\theta} - \frac{1}{\gamma}-\varepsilon)\theta} \Biggl(\sum\limits_{s_{m}=\nu_{m}}^{\infty}...\sum\limits_{s_{1}=\nu_{1}}^{\infty}\|\delta_{\overline{s}}(f)\|_{p, \tau}^{\gamma}
\Biggr)^{\frac{\theta}{\gamma}}\Biggr)^{\frac{1}{\theta}}
\\
=C \Biggl(\sum\limits_{\nu_{m} =1}^{\infty}...\sum\limits_{\nu_{1} =1}^{\infty}\prod_{j=1}^{m}\nu_{j}^{b_{j}\theta} \Biggl(\sum\limits_{s_{m}=\nu_{m}}^{\infty}...\sum\limits_{s_{1}=\nu_{1}}^{\infty}\|\delta_{\overline{s}}(f)\|_{p, \tau}^{\gamma}
\Biggr)^{\frac{\theta}{\gamma}}\Biggr)^{\frac{1}{\theta}},
\end{multline}
in the case  $0<\theta< \gamma=\max\{2, p\}$ and $b_{j}> \frac{1}{\theta} - \frac{1}{\gamma}$, $j=1, \ldots , m$. Further, from inequalities (2. 23) and (2. 31) we obtain \begin{multline}\label{eq2 32}
\Biggl(\sum\limits_{s_{m} =1}^{\infty}...\sum\limits_{s_{1} =1}^{\infty}\prod_{j=1}^{m}s_{j}^{(b_{j} + \frac{1}{\gamma})\theta}\|\delta_{\overline{s}}(f)\|_{p, \tau}^{\theta}\Biggr)^{\frac{1}{\theta}}
\\
\leq C \Biggl(\sum\limits_{\nu_{m} =1}^{\infty}...\sum\limits_{\nu_{1} =1}^{\infty}\prod_{j=1}^{m}\nu_{j}^{b_{j}\theta} \Biggl\|\Biggl(\sum\limits_{s_{m}=\nu_{m}}^{\infty}...\sum\limits_{s_{1}=\nu_{1}}^{\infty}|\delta_{\overline{s}}(f)|^{2}\Biggr)^{\frac{1}{2}}\Biggr\|_{p, \tau}^{\theta}\Biggr)^{\frac{1}{\theta}},
\end{multline}
in the case  $0<\theta< \gamma$ and $b_{j}> \frac{1}{\theta} - \frac{1}{\gamma}$, $j=1, \ldots , m$. 
 Now,  from inequality (2. 32) and Theorem 2, it follows that
 $S_{p, \tau, \theta}^{0, \overline{b}}\mathbf{B} \subset S_{p, \tau, \theta}^{0, \overline{u}}B$, in the case   $0<\theta< \gamma$ and $b_{j}> \frac{1}{\theta} - \frac{1}{\gamma}$, $j=1, \ldots , m$.
 \hfill $\Box$   

\begin{theorem}\label{th5}   
  Let $1 < p <+ \infty$, $1< \tau_{2}\leq \tau_{1} < \infty$, $0< \theta\leq \infty$ and the numbers $b_{j}^{(i)}> - \frac{1}{\theta_{i}}$, for $j=1,\ldots, m$, $\overline{b}^{(i)}=(b_{1}^{(i)}, \ldots , b_{m}^{(i)})$, $i=1, 2$. Then
    
 1.    $S_{p, \tau_{2}, \theta}^{0, \overline{b}^{(2)}}B \subset S_{p, \tau_{1}, \theta}^{0, \overline{b}^{(1)}}B$, if   $b_{j}^{(1)}< b_{j}^{(2)}$, $j=1,\ldots, m$;
 
 2. if  $1< \tau_{2}< \tau_{1}  < \infty$, $0< \theta_{2}< \theta_{1}\leq \infty$ and
\begin{equation}\label{eq2 33} 
\Biggl(\sum\limits_{s_{m} =1}^{\infty}...\sum\limits_{s_{1} =1}^{\infty}\prod_{j=1}^{m}s_{j}^{(b_{j}^{(2)}-b_{j}^{(1)})\theta_{2}\eta^{'}}\Bigl(\sum\limits_{j =1}^{m}(s_{j} + 1)\Bigr)^{(\frac{1}{\tau_{2}}-\frac{1}{\tau_{1}})\theta_{2}\eta^{'}}\Biggr)^{\frac{1}{\theta_{2}\eta^{'}}} < \infty, 
\end{equation}
then 
 $S_{p, \tau_{1}, \theta_{1}}^{0, \overline{b}^{(1)}}B \subset S_{p, \tau_{2}, \theta_{2}}^{0, \overline{b}^{(2)}}B$, where 
  $\eta = \frac{\theta_{1}}{\theta_{2}}\, \, \,\,$, $\frac{1}{\eta} + \frac{1}{\eta^{`}} = 1$;

3. if  
$b_{j}^{(1)}+\frac{1}{\tau_{1}}  >b_{j}^{(2)}+\frac{1}{\tau_{2}}$, for $j=1,\ldots, m$, $1< \tau_{2}< \tau_{1}  < \infty$, $0< \theta_{2}< \theta_{1}\leq \infty$, 
\begin{equation*} 
\sum\limits_{l_{m} =0}^{\infty}...\sum\limits_{l_{1} =0}^{\infty}\prod_{j=1}^{m} 2^{l_{j}(b_{j}^{(2)} - b_{j}^{(1)} - \frac{1}{\theta_{1}} + \frac{1}{\theta_{2}})\theta_{2}\eta^{`}}\Bigl(\sum_{j=1}^{m}2^{l_{j}}\Bigr)^{(\frac{1}{\tau_{2}} - \frac{1}{\tau_{1}})\theta_{2}\eta^{`}} < \infty,
\end{equation*}
then    $S_{p, \tau_{1}, \theta_{1}}^{0, \overline{b}^{(1)}}\mathbf{B} \subset S_{p, \tau_{2}, \theta_{2}}^{0, \overline{b}^{(2)}}\mathbf{B}$.
\end{theorem}

\proof It is known that if $1< \tau_{2}\leq\tau_{1} <\infty$, then $\|f\|_{p, \tau_{1}}\ll\|f\|_{p, \tau_{2}}$, for $f\in L_{p, \tau_{2}}$, $1 < p <+ \infty$ \cite[theorem 3. 11]{1}.
The first statement of the theorem follows from this inequality.

We will prove the second point. Let $f\in S_{p, \tau_{1}, \theta_{1}}^{0, \overline{b}^{(1)}}B$.
Since
$1< \tau_{2}< \tau_{1} < \infty$, then according to the inequality of different metrics for a trigonometric polynomial in the Lorentz space \cite[lemma 3.1]{27} we have
\begin{multline*}
\Biggl(\sum\limits_{s_{m} =1}^{\infty}...\sum\limits_{s_{1} =1}^{\infty}\prod_{j=1}^{m}s_{j}^{b_{j}^{(2)}\theta_{2}}\|\delta_{\overline{s}}(f)\|_{p, \tau_{2}}^{\theta_{2}}\Biggr)^{\frac{1}{\theta_{2}}}
\\
\ll \Biggl(\sum\limits_{s_{m} =1}^{\infty}...\sum\limits_{s_{1} =1}^{\infty}\prod_{j=1}^{m}s_{j}^{b_{j}^{(2)}\theta_{2}}\Bigl(\sum\limits_{j =1}^{m}(s_{j} + 1)\Bigr)^{(\frac{1}{\tau_{2}}-\frac{1}{\tau_{1}})\theta_{2}}\|\delta_{\overline{s}}(f)\|_{p, \tau_{2}}^{\theta_{2}}\Biggr)^{\frac{1}{\theta_{2}}}.
\end{multline*}
Since $0< \theta_{2}< \theta_{1}$, then applying the H\"{o}lder inequality for $\eta=\frac{\theta_{1}}{\theta_{2}}$, $\, \, \frac{1}{\eta}+\frac{1}{\beta^{'}}=1$, hence we get
 \begin{multline*}
\Biggl(\sum\limits_{s_{m} =1}^{\infty}...\sum\limits_{s_{1} =1}^{\infty}\prod_{j=1}^{m}s_{j}^{b_{j}^{(2)}\theta_{2}}\|\delta_{\overline{s}}(f)\|_{p, \tau_{2}}^{\theta_{2}}\Biggr)^{\frac{1}{\theta_{2}}}\ll \Biggl(\sum\limits_{s_{m} =1}^{\infty}...\sum\limits_{s_{1} =1}^{\infty}\prod_{j=1}^{m}s_{j}^{b_{j}^{(1)}\theta_{1}}\|\delta_{\overline{s}}(f)\|_{p, \tau_{1}}^{\theta_{1}}\Biggr)^{\frac{1}{\theta_{1}}} 
\\
\times
\Biggl(\sum\limits_{s_{m} =1}^{\infty}...\sum\limits_{s_{1} =1}^{\infty}\prod_{j=1}^{m}s_{j}^{(b_{j}^{(2)}-b_{j}^{(1)})\theta_{2}\eta^{'}}\Bigl(\sum\limits_{j =1}^{m}(s_{j} + 1)\Bigr)^{(\frac{1}{\tau_{2}}-\frac{1}{\tau_{1}})\theta_{2}\eta^{'}} \Biggr)^{\frac{1}{\theta_{2}\eta^{'}}}.
\end{multline*}
According to condition (2. 33) and the definition of the space $S_{p, \tau_{1}, \theta_{1}}^{0, \overline{b}^{(1)}}B$, from here we get $S_{p, \tau_{1}, \theta_{1}}^{0, \overline{b}^{(1)}}B \subset S_{p, \tau_{2}, \theta_{2}}^{0, \overline{b}^{(2)}}B$.

Statement 3 is proved similarly, based on Theorem 3.
\hfill $\Box$  

\begin{rem}
  From Theorem 4 and the first statement of Theorem 5, it follows that the spaces $S_{p, \tau, \theta}^{0, \overline{b}}\mathbf{B}$ and $S_{p, \tau, \theta}^{0, \overline{b}}B$ does not co-occur with $b_{j}> - \frac{1}{\theta}$, for $j=1,\ldots, m$.
\end{rem} 

\section*{Conclusion}
The analog of the space $\mathbf{B}_{p, \theta}^{0, b}(\mathbb{T}^{m})$ in the Lorentz-Zygmund space is defined and investigated in \cite{21}, and the spaces $B_{p, \theta}^{0,b}(\mathbb{T}^{m})$ in the generalized Lorentz space in \cite{28}.

\end{document}